\documentclass[11pt]{article}
\UseRawInputEncoding
\usepackage{amscd}
\usepackage{amssymb,amsfonts,amsmath,psfrag,amscd,cite}
\usepackage{amsmath}
  \usepackage{paralist}
  \usepackage{graphics}
  \usepackage{epsfig}
  \usepackage[colorlinks=true]{hyperref}
  \usepackage[all]{xy}
\usepackage{graphicx}
\newtheorem{theo}{Theorem}
\newtheorem{prop}[theo]{Proposition}

\newtheorem{defi}[theo]{Definition}
\newtheorem{lemm}[theo]{Lemma}
\newtheorem{rema}[theo]{Remark}

\makeatletter \@addtoreset{equation}{section}
\@addtoreset{theo}{section} \makeatother

\begin{document}
\date{}
\title{Hamilton-Jacobi Equations for Nonholonomic Magnetic Hamiltonian Systems }
\author{Hong Wang \\
School of Mathematical Sciences and LPMC,\\
Nankai University,  Tianjin 300071, P.R.China\\
E-mail: hongwang@nankai.edu.cn\\\\
\emph{In Memory of Great Geometer Shiing Shen Chern} \\
June 15, 2022} \maketitle

{\bf Abstract:} In order to describe the impact of different geometric structures
and constraints for the dynamics of a Hamiltonian system,
in this paper, for a magnetic Hamiltonian system
defined by a magnetic symplectic form, we first
drive precisely the geometric constraint conditions of
magnetic symplectic form for the magnetic Hamiltonian vector field.
which are called the Type I and Type II of Hamilton-Jacobi equation.
Secondly, for the magnetic Hamiltonian system with nonholonomic constraint,
we first define a distributional magnetic Hamiltonian system,
then derive its two types of Hamilton-Jacobi equation.
Moreover, we generalize the above results
to nonholonomic reducible magnetic Hamiltonian system with symmetry.
We define a nonholonomic reduced distributional
magnetic Hamiltonian system, and prove two types of
Hamilton-Jacobi theorem. These research work reveal the deeply internal
relationships of the magnetic symplectic structure, nonholonomic constraint,
the distributional two-form, and the dynamical vector field of
the nonholonomic magnetic Hamiltonian system.\\



\tableofcontents

\section{Introduction}

It is well-known that Hamilton-Jacobi theory is an important research subject
in mathematics and analytical mechanics,
see Abraham and Marsden \cite{abma78}, Arnold
\cite{ar89} and Marsden and Ratiu \cite{mara99},
and the Hamilton-Jacobi equation is
also fundamental in the study of the quantum-classical relationship
in quantization, and it also plays an important role
in the study of stochastic dynamical systems, see
Woodhouse \cite{wo92}, Ge and Marsden \cite{gema88},
and L\'{a}zaro-Cam\'{i} and Ortega \cite{laor09}.
For these reasons it is described as a useful tool in the study of
Hamiltonian system theory, and has been extensively developed in
past many years and become one of the most active subjects in the
study of modern applied mathematics and analytical mechanics.\\

Just as we have known that Hamilton-Jacobi theory from the
variational point of view is originally developed by Jacobi in 1866,
which state that the integral of Lagrangian of a mechanical system along the
solution of its Euler-Lagrange equation satisfies the
Hamilton-Jacobi equation. The classical description of this problem
from the generating function and the geometrical point of view is
given by Abraham and Marsden in \cite{abma78} as follows:
Let $Q$ be a smooth manifold and $TQ$
the tangent bundle, $T^* Q$ the cotangent bundle with a canonical
symplectic form $\omega$ and the projection $\pi_Q: T^* Q
\rightarrow Q $ induces the map $T\pi_{Q}: TT^* Q \rightarrow TQ. $
\begin{theo}
Assume that the triple $(T^*Q,\omega,H)$ is a Hamiltonian system
with Hamiltonian vector field $X_H$, and $W: Q\rightarrow
\mathbb{R}$ is a given generating function. Then the following two assertions
are equivalent:\\
\noindent $(\mathrm{i})$ For every curve $\sigma: \mathbb{R}
\rightarrow Q $ satisfying $\dot{\sigma}(t)= T\pi_Q
(X_H(\mathbf{d}W(\sigma(t))))$, $\forall t\in \mathbb{R}$, then
$\mathbf{d}W \cdot \sigma $ is an integral curve of the Hamiltonian
vector field $X_H$.\\
\noindent $(\mathrm{ii})$ $W$ satisfies the Hamilton-Jacobi equation
$H(q^i,\frac{\partial W}{\partial q^i})=E, $ where $E$ is a
constant.
\end{theo}

From the proof of the above theorem given in
Abraham and Marsden \cite{abma78}, we know that
the assertion $(\mathrm{i})$ with equivalent to
Hamilton-Jacobi equation $(\mathrm{ii})$ by the generating function,
gives a geometric constraint condition of the canonical symplectic form
on the cotangent bundle $T^*Q$
for Hamiltonian vector field of the system.
Thus, the Hamilton-Jacobi equation reveals the deeply internal relationships of
the generating function, the canonical symplectic form
and the dynamical vector field of a Hamiltonian system.\\

Now, it is a natural problem how to generalize Theorem 1.1 to fit
the nonholonomic systems and their reduced systems.
Note that if take that $\gamma=\mathbf{d}W$ in
the above Theorem 1.1, then $\gamma$ is a closed one-form on $Q$, and
the equation $\mathbf{d}(H \cdot \mathbf{d}W)=0$ is equivalent to
the Hamilton-Jacobi equation $H(q^i,\frac{\partial W}{\partial
q^i})=E$, where $E$ is a constant, which is called the classical
Hamilton-Jacobi equation. This result is used the
formulation of a geometric version of Hamilton-Jacobi theorem for
Hamiltonian system, see Cari\~{n}ena et al \cite{cagrmamamuro06,
cagrmamamuro10}. Moreover, note that Theorem 1.1 is also
generalized in the context of time-dependent Hamiltonian system by
Marsden and Ratiu in \cite{mara99}, and the Hamilton-Jacobi equation may
be regarded as a nonlinear partial differential equation for some
generating function $S$. Thus, the problem is become how to choose a
time-dependent canonical transformation $\Psi: T^*Q\times \mathbb{R}
\rightarrow T^*Q\times \mathbb{R}, $ which transforms the dynamical
vector field of a time-dependent Hamiltonian system to equilibrium,
such that the generating function $S$ of $\Psi$ satisfies the
time-dependent Hamilton-Jacobi equation. In particular, for the
time-independent Hamiltonian system, ones may look for a symplectic
map as the canonical transformation. This work offers an important
idea that one can use the dynamical vector field of a Hamiltonian
system to describe Hamilton-Jacobi equation. In consequence, if assume that
$\gamma: Q \rightarrow T^*Q$ is a closed one-form on $Q$, and define
that $X_H^\gamma = T\pi_{Q}\cdot X_H \cdot \gamma$, where $X_{H}$ is
the dynamical vector field of Hamiltonian system $(T^*Q,\omega,H)$,
then the fact that $X_H^\gamma$ and $X_H$ are $\gamma$-related, that
is, $T\gamma\cdot X_H^\gamma= X_H\cdot \gamma$ is equivalent that
$\mathbf{d}(H \cdot \gamma)=0, $ which is given in Cari\~{n}ena et
al \cite{cagrmamamuro06, cagrmamamuro10}.
Motivated by the above research work, Wang in \cite{wa17} prove
an important lemma, which is
a modification for the corresponding result of Abraham and Marsden
in \cite{abma78}, such that we can derive precisely
the geometric constraint conditions of
the regular reduced symplectic forms for the
dynamical vector fields of a regular reducible Hamiltonian
system on the cotangent bundle of a configuration manifold,
which are called the Type I and Type II of Hamilton-Jacobi equation,
because they are the development of the above classical Hamilton-Jacobi equation
given by Theorem 1.1, see Abraham and Marsden \cite{abma78}
and Wang \cite{wa17}. Moreover, Le\'{o}n and Wang in \cite{lewa15}
generalize the above results to
the nonholonomic Hamiltonian system and
the nonholonomic reducible Hamiltonian system
on a cotangent bundle, by using the distributional Hamiltonian system
and the reduced distributional Hamiltonian system.\\

In order to describe the impact of different geometric structures
and constraints for the dynamics of a Hamiltonian system,
in the following we first consider the magnetic Hamiltonian system.
Define a magnetic symplectic form $\omega^B= \omega-
\pi^*_Q B,$  and the $\pi_Q^*B$ is called a magnetic term on $T^*Q$,
where $\omega$ is the usual canonical symplectic
form on $T^*Q$, and $B$ is the closed two-form on $Q$,
and the map $\pi^*_Q: T^* Q \rightarrow T^*T^*Q $.
A magnetic Hamiltonian system is a Hamiltonian system defined by the
magnetic symplectic form, which is a canonical Hamiltonian system
coupling the action of a magnetic field $B$. Under the impact of
magnetic term $\pi_Q^*B$, the magnetic symplectic form $\omega^B$,
in general, is not the canonical symplectic form on $T^*Q$,
we cannot prove the Hamilton-Jacobi theorem for the magnetic Hamiltonian system
just like same as the above Theorem 1.1. We have to look for a new way.
In this paper, we first drive precisely the geometric constraint conditions of
magnetic symplectic form for the magnetic Hamiltonian vector field.
These conditions are called the Type I and Type II of Hamilton-Jacobi equation,
which are the development of the Type I and Type II of
Hamilton-Jacobi equation for a Hamiltonian system
given in Wang \cite{wa17}.\\

Next, we consider the magnetic Hamiltonian system
with nonholonomic constraint, which is called the nonholonomic
magnetic Hamiltonian system. In mechanics,
it is very often that many systems have constraints,
and usually, under the restriction given by nonholonomic constraint,
in general, the dynamical vector field of a nonholonomic magnetic Hamiltonian system
may not be Hamiltonian. Thus, we can not
describe the Hamilton-Jacobi equations for nonholonomic magnetic
Hamiltonian system from the viewpoint of generating function
as in the classical Hamiltonian case, that is,
we cannot prove the Hamilton-Jacobi
theorem for the nonholonomic magnetic Hamiltonian system,
just like same as the above Theorem 1.1. In this paper,
by analyzing carefully the structure for the nonholonomic dynamical
vector field, we first give a geometric formulation of the distributional
magnetic Hamiltonian system for the nonholonomic magnetic Hamiltonian system,
which is determined by a non-degenerate distributional two-form induced
from the magnetic symplectic form. The distributional
magnetic Hamiltonian system is not Hamiltonian,
however, it is a dynamical system closely related to a
magnetic Hamiltonian system. Then we drive precisely
two types of Hamilton-Jacobi equation for the distributional magnetic
Hamiltonian system, which are the development of the Type I and Type II of
Hamilton-Jacobi equation for a distributional Hamiltonian system
given in Le\'{o}n and Wang \cite{lewa15}.\\

Thirdly, it is a natural problem to consider the nonholonomic magnetic Hamiltonian system
with symmetry. In this paper, we generalize the above results
to nonholonomic reducible magnetic Hamiltonian system with symmetry.
By using the method of nonholonomic reduction
given in Bates and $\acute{S}$niatycki \cite{basn93},
and analyzing carefully the structure for the nonholonomic reduced
dynamical vector field, we first give a geometric formulation of
the nonholonomic reduced distributional magnetic Hamiltonian system.
Since the nonholonomic reduced distributional
magnetic Hamiltonian system is not yet Hamiltonian,
but, it is a dynamical system closely related to a
magnetic Hamiltonian system.
Then we can derive precisely the geometric constraint conditions of
the non-degenerate, and nonholonomic reduced distributional two-form
for the nonholonomic reducible dynamical vector field,
that is, the two types of Hamilton-Jacobi equation for the
nonholonomic reduced distributional magnetic Hamiltonian system.
These research work reveal the deeply internal
relationships of the magnetic symplectic structure, nonholonomic constraint,
the induced (resp. reduced) distributional two-forms,
and the dynamical vector fields of the nonholonomic magnetic
Hamiltonian system.\\

The paper is organized as follows. In section 2 we first
give some definitions and basic facts
about the magnetic Hamiltonian system, the nonholonomic
constraint, the nonholonomic magnetic Hamiltonian system
and the distributional magnetic Hamiltonian system,
which will be used in subsequent sections.
In section 3, for a magnetic Hamiltonian system
defined by a magnetic symplectic form, we first
drive precisely the geometric constraint conditions of
magnetic symplectic form for the magnetic Hamiltonian vector field.
which are called the Type I and Type II of Hamilton-Jacobi equation.
In section 4, we derive two types of Hamilton-Jacobi equation for a
distributional magnetic Hamiltonian system, by the analysis and calculation
in detail. The nonholonomic reducible magnetic Hamiltonian system with
symmetry is considered in section 5, and
derive precisely the geometric constraint conditions of
the non-degenerate, and nonholonomic reduced distributional two-form
for the nonholonomic reducible dynamical vector field,
that is, the two types of Hamilton-Jacobi
equations for the nonholonomic reduced
 distributional magnetic Hamiltonian system.
These research work develop Hamilton-Jacobi theory
for the nonholonomic magnetic Hamiltonian system, as well as with symmetry,
and make us have much deeper understanding
and recognition for the structures of the nonholonomic
magnetic Hamiltonian systems.

\section{Nonholonomic Magnetic Hamiltonian System}

In this section we first give some definitions and basic facts
about the magnetic Hamiltonian system, the nonholonomic
constraint and the nonholonomic magnetic Hamiltonian system.
Moreover, by analyzing carefully the structure for the nonholonomic dynamical
vector field, we give a geometric formulation of distributional magnetic
Hamiltonian system, which is
determined by a non-degenerate distributional two-form induced
from the magnetic symplectic form.
All of them will be used in subsequent sections.\\

Let $Q$ be an $n$-dimensional smooth manifold and $TQ$
the tangent bundle, $T^* Q$ the cotangent bundle with a canonical
symplectic form $\omega$ and the projection $\pi_Q: T^* Q
\rightarrow Q $ induces the map $\pi^*_{Q}: T^* Q \rightarrow T^*T^*Q. $
We consider the magnetic symplectic form
$\omega^B= \omega- \pi^*_Q B,$ where $\omega$ is the canonical symplectic
form on $T^*Q$, and $B$ is the closed two-form on $Q$,
and the $\pi_Q^*B$ is called a magnetic term on $T^*Q$.
A magnetic Hamiltonian system is a triple $(T^\ast Q,\omega^B,H)$,
which is a Hamiltonian system defined by the
magnetic symplectic form $\omega^B$, that is,
a canonical Hamiltonian system
coupling the action of a magnetic field $B$. For a given Hamiltonian $H$,
the dynamical vector field $X^B_H$, which is called
the magnetic Hamiltonian vector field,
satisfies the magnetic Hamilton's equation, that is,
$\mathbf{i}_{X^B_{H} }\omega^B= \mathbf{d}H $.
In canonical cotangent bundle coordinates, for any $q \in Q
, \; (q,p)\in T^* Q, $ we have that
$$
\omega=\sum^n_{i=1} \mathbf{d}q^i \wedge \mathbf{d}p_i ,
\;\;\;\;\;\; B=\sum^n_{i,j=1}B_{ij}\mathbf{d}q^i \wedge
\mathbf{d}q^j ,\;\;\; \mathbf{d}B=0, $$
$$\omega^B= \omega -\pi_Q^*B=\sum^n_{i=1} \mathbf{d}q^i \wedge
\mathbf{d}p_i- \sum^n_{i,j=1}B_{ij}\mathbf{d}q^i \wedge
\mathbf{d}q^j,
$$
and the magnetic Hamiltonian vector field $X^B_H$ with respect to
the magnetic symplectic form $\omega^B$ can be expressed that
$$
X^B_H= \sum^n_{i=1} (\frac{\partial H}{\partial
p_i}\frac{\partial}{\partial q^i} - \frac{\partial H}{\partial
q^i}\frac{\partial}{\partial p_i})-
\sum^n_{i,j=1}B_{ij}\frac{\partial H}{\partial
p_j}\frac{\partial}{\partial p_i}.
$$
See Marsden et al. \cite{mamiorpera07}.\\

In order to describe the nonholonomic magnetic Hamiltonian system,
in the following we first give the completeness and regularity
conditions for nonholonomic constraints of a mechanical system,
see Le\'{o}n and Wang \cite{lewa15}. In fact,
in order to describe the dynamics of a nonholonomic mechanical system,
we need some restriction conditions for nonholonomic constraints of
the system. At first, we note that the set of Hamiltonian vector fields
forms a Lie algebra with respect to the Lie bracket, since
$X_{\{f,g\}}=-[X_f, X_g]. $ But, the Lie bracket operator, in
general, may not be closed on the restriction of a nonholonomic
constraint. Thus, we have to give the following completeness
condition for nonholonomic constraints of a system.\\

{\bf $\mathcal{D}$-completeness } Let $Q$ be a smooth manifold and
$TQ$ its tangent bundle. A distribution $\mathcal{D} \subset TQ$ is
said to be {\bf completely nonholonomic} (or bracket-generating) if
$\mathcal{D}$ along with all of its iterated Lie brackets
$[\mathcal{D},\mathcal{D}], [\mathcal{D}, [\mathcal{D},\mathcal{D}]],
\cdots ,$ spans the tangent bundle $TQ$. Moreover, we consider a
nonholonomic mechanical system on $Q$, which is
given by a Lagrangian function $L: TQ \rightarrow \mathbb{R}$
subject to constraints determined by a nonholonomic
distribution $\mathcal{D}\subset TQ$ on the configuration manifold $Q$.
Then the nonholonomic system is said to be {\bf completely nonholonomic},
if the distribution $\mathcal{D} \subset TQ$ determined by the nonholonomic
constraints is completely nonholonomic.\\

{\bf $\mathcal{D}$-regularity } In the following we always assume
that $Q$ is a smooth manifold with coordinates $(q^i)$, and $TQ$ its
tangent bundle with coordinates $(q^i,\dot{q}^i)$, and $T^\ast Q$
its cotangent bundle with coordinates $(q^i,p_j)$, which are the
canonical cotangent coordinates of $T^\ast Q$ and $\omega=
dq^{i}\wedge dp_{i}$ is canonical symplectic form on $T^{\ast}Q$. If
the Lagrangian $L: TQ \rightarrow \mathbb{R}$ is hyperregular, that
is, the Hessian matrix
$(\partial^2L/\partial\dot{q}^i\partial\dot{q}^j)$ is nondegenerate
everywhere, then the Legendre transformation $FL: TQ \rightarrow T^*
Q$ is a diffeomorphism. In this case the Hamiltonian $H: T^* Q
\rightarrow \mathbb{R}$ is given by $H(q,p)=\dot{q}\cdot
p-L(q,\dot{q}) $ with Hamiltonian vector field $X_H$,
which is defined by the Hamilton's equation
$\mathbf{i}_{X_H}\omega=\mathbf{d}H$, and
$\mathcal{M}=\mathcal{F}L(\mathcal{D})$ is a constraint submanifold
in $T^* Q$. In particular, for the nonholonomic constraint
$\mathcal{D}\subset TQ$, the Lagrangian $L$ is said to be {\bf
$\mathcal{D}$-regular}, if the restriction of Hessian matrix
$(\partial^2L/\partial\dot{q}^i\partial\dot{q}^j)$ on $\mathcal{D}$
is nondegenerate everywhere. Moreover, a nonholonomic system is said
to be {\bf $\mathcal{D}$-regular}, if its Lagrangian $L$ is {\bf
$\mathcal{D}$-regular}. Note that the restriction of a positive
definite symmetric bilinear form to a subspace is also positive
definite, and hence nondegenerate. Thus, for a simple nonholonomic
mechanical system, that is, whose Lagrangian is the total kinetic
energy minus potential energy, it is {\bf $\mathcal{D}$-regular }
automatically.\\

A nonholonomic magnetic Hamiltonian system is a 4-tuple
$(T^\ast Q,\omega^B,\mathcal{D},H)$,
which is a magnetic Hamiltonian system with a
$\mathcal{D}$-completely and $\mathcal{D}$-regularly nonholonomic
constraint $\mathcal{D} \subset TQ$. Under
the restriction given by constraint, in general, the dynamical
vector field of a nonholonomic magnetic Hamiltonian system may not be
magnetic Hamiltonian, however the system is a dynamical system
closely related to a magnetic Hamiltonian system.
In the following we shall derive a distributional magnetic
Hamiltonian system of the nonholonomic
magnetic Hamiltonian system $(T^*Q,\omega^B,\mathcal{D},H)$,
by analyzing carefully the structure for the nonholonomic
dynamical vector field similar to the method in Le\'{o}n and Wang \cite{lewa15}.
It is worthy of noting that the leading distributional Hamiltonian system
is also called a semi-Hamiltonian system in Patrick \cite{pa07}.\\

We consider that the constraint submanifold
$\mathcal{M}=\mathcal{F}L(\mathcal{D})\subset T^*Q$ and
$i_{\mathcal{M}}: \mathcal{M}\rightarrow T^*Q $ is the inclusion,
the symplectic form $\omega^B_{\mathcal{M}}= i_{\mathcal{M}}^* \omega^B $,
is induced from the magnetic symplectic form $\omega^B$ on $T^* Q$.
We define the distribution $\mathcal{F}$ as the pre-image of the nonholonomic
constraints $\mathcal{D}$ for the map $T\pi_Q: TT^* Q \rightarrow TQ$,
that is, $\mathcal{F}=(T\pi_Q)^{-1}(\mathcal{D})\subset TT^*Q,
$ which is a distribution along $\mathcal{M}$, and
$\mathcal{F}^\circ:=\{\alpha \in T^*T^*Q | <\alpha,v>=0, \; \forall
v\in TT^*Q \}$ is the annihilator of $\mathcal{F}$ in
$T^*T^*Q_{|\mathcal{M}}$. We consider the following nonholonomic
constraints condition
\begin{align} (\mathbf{i}_X \omega^B -\mathbf{d}H) \in \mathcal{F}^\circ,
\;\;\;\;\;\; X \in T \mathcal{M},
\label{2.1} \end{align} from Cantrijn et al.
\cite{calemama99}, we know that there exists an unique nonholonomic
vector field $X_n$ satisfying the above condition $(2.1)$, if the
admissibility condition $\mathrm{dim}\mathcal{M}=
\mathrm{rank}\mathcal{F}$ and the compatibility condition
$T\mathcal{M}\cap \mathcal{F}^\bot= \{0\}$ hold, where
$\mathcal{F}^\bot$ denotes the magnetic symplectic orthogonal of
$\mathcal{F}$ with respect to the magnetic symplectic form
$\omega^B$ on $T^*Q$. In particular, when we consider the Whitney sum
decomposition $T(T^*Q)_{|\mathcal{M}}=T\mathcal{M}\oplus
\mathcal{F}^\bot$ and the canonical projection $P:
T(T^*Q)_{|\mathcal{M}} \rightarrow T\mathcal{M}$,
then we have that $X_n= P(X^B_H)$.\\

From the condition (2.1) we know that the nonholonomic vector field,
in general, may not be magnetic Hamiltonian, because of the restriction
of nonholonomic constraints. But, we hope to study the dynamical
vector field of nonholonomic magnetic Hamiltonian system by using the similar
method of studying magnetic Hamiltonian vector field.
From Le\'{o}n and Wang \cite{lewa15} and Bates
and $\acute{S}$niatycki \cite{basn93}, by using the
similar method, we can define the
distribution $ \mathcal{K}=\mathcal {F}\cap T\mathcal{M}.$ and
$\mathcal{K}^\bot=\mathcal {F}^\bot\cap T\mathcal{M}, $ where
$\mathcal{K}^\bot$ denotes the magnetic symplectic orthogonal of
$\mathcal{K}$ with respect to the magnetic symplectic form
$\omega^B$, and the admissibility condition $\mathrm{dim}\mathcal{M}=
\mathrm{rank}\mathcal{F}$ and the compatibility condition
$T\mathcal{M}\cap \mathcal{F}^\bot= \{0\}$ hold, then we know that the
restriction of the symplectic form $\omega^B_{\mathcal{M}}$ on
$T^*\mathcal{M}$ fibrewise to the distribution $\mathcal{K}$, that
is, $\omega^B_\mathcal{K}= \tau_{\mathcal{K}}\cdot
\omega^B_{\mathcal{M}}$ is non-degenerate, where $\tau_{\mathcal{K}}$
is the restriction map to distribution $\mathcal{K}$. It is worthy
of noting that $\omega^B_\mathcal{K}$ is not a true two-form on a
manifold, so it does not make sense to speak about it being closed.
We call $\omega^B_\mathcal{K}$ as a distributional two-form to avoid
any confusion. Because $\omega^B_\mathcal{K}$ is non-degenerate as a
bilinear form on each fibre of $\mathcal{K}$, there exists a vector
field $X^B_{\mathcal{K}}$ on $\mathcal{M}$ which takes values in the
constraint distribution $\mathcal{K}$,
such that the distributional magnetic Hamiltonian equation
\begin{align}
\mathbf{i}_{X^B_\mathcal{K}}\omega^B_{\mathcal{K}}
=\mathbf{d}H_\mathcal{K}
\label{2.2} \end{align}
holds, where $\mathbf{d}H_\mathcal{K}$ is the restriction of
$\mathbf{d}H_\mathcal{M}$ to $\mathcal{K}$,
and the function $H_{\mathcal{K}}$ satisfies
$\mathbf{d}H_{\mathcal{K}}= \tau_{\mathcal{K}}\cdot \mathbf{d}H_{\mathcal {M}}$,
and $H_\mathcal{M}= \tau_{\mathcal{M}}\cdot H$ is the restriction of $H$ to
$\mathcal{M}$. Moreover, from the distributional magnetic Hamiltonian equation (2.2),
we have that $X^B_{\mathcal{K}}= \tau_{\mathcal{K}}\cdot X^B_H.$
Then the triple $(\mathcal{K},\omega^B_{\mathcal{K}},H_{\mathcal{K}})$
is a distributional magnetic Hamiltonian system of the nonholonomic
magnetic Hamiltonian system $(T^*Q,\omega^B,\mathcal{D},H)$.
Thus, the geometric formulation of the distributional magnetic
Hamiltonian system may be summarized as follows.

\begin{defi} (Distributional Magnetic Hamiltonian System)
Assume that the 4-tuple $(T^*Q,\omega^B,\\ \mathcal{D},H)$ is a
$\mathcal{D}$-completely and $\mathcal{D}$-regularly nonholonomic
magnetic Hamiltonian system, where the
magnetic symplectic form $\omega^B= \omega- \pi_Q^*B $ on $T^*Q$,
and $\omega$ is the canonical symplectic form on $T^* Q$
and $B$ is a closed two-form on $Q$, and
$\mathcal{D}\subset TQ$ is a
$\mathcal{D}$-completely and $\mathcal{D}$-regularly nonholonomic
constraint of the system. If there exist a distribution
$\mathcal{K}$, an associated non-degenerate distributional two-form
$\omega^B_{\mathcal{K}}$ induced by the magnetic symplectic form $\omega^B$
and a vector field $X^B_\mathcal {K}$ on the
constraint submanifold $\mathcal{M}=\mathcal{F}L(\mathcal{D})\subset
T^*Q$, such that the distributional magnetic Hamiltonian equation
$\mathbf{i}_{X^B_\mathcal{K}}\omega^B_{\mathcal{K}}=\mathbf{d}H_\mathcal
{K}$ holds, where $\mathbf{d}H_\mathcal{K}$ is the restriction of
$\mathbf{d}H_\mathcal{M}$ to $\mathcal{K}$
and the function $H_{\mathcal{K}}$ satisfies
$\mathbf{d}H_{\mathcal{K}}= \tau_{\mathcal{K}}\cdot \mathbf{d}H_{\mathcal {M}}$
as defined above,
then the triple $(\mathcal{K},\omega^B_{\mathcal{K}},H_{\mathcal{K}})$
is called a distributional magnetic Hamiltonian system of the nonholonomic
magnetic Hamiltonian system $(T^*Q,\omega^B,\mathcal{D},H)$, and $X^B_\mathcal
{K}$ is called a nonholonomic dynamical
vector field of the distributional magnetic Hamiltonian system
$(\mathcal{K},\omega^B_{\mathcal {K}},H_{\mathcal{K}})$. Under the above circumstances, we refer to
$(T^*Q,\omega^B,\mathcal{D},H)$ as a nonholonomic magnetic Hamiltonian system
with an associated distributional magnetic Hamiltonian system
$(\mathcal{K},\omega^B_{\mathcal {K}},H_{\mathcal{K}})$.
\end{defi}

Moreover, in section 5, we consider the nonholonomic magnetic Hamiltonian system
with symmetry. By using the similar method for nonholonomic reduction
given in Bates and $\acute{S}$niatycki \cite{basn93}
and Le\'{o}n and Wang \cite{lewa15},
and analyzing carefully the structure for the nonholonomic reduced
dynamical vector field, we also give a geometric formulation of
the nonholonomic reduced distributional magnetic Hamiltonian system.

\section{Hamilton-Jacobi Equation of Magnetic Hamiltonian System}

In order to describe the impact of different geometric structures and constraints
for the Hamilton-Jacobi theory, in this paper, we shall give two
types of Hamilton-Jacobi equations for the magnetic Hamiltonian system,
the distributional magnetic Hamiltonian system and the nonholonomic
reduced distributional magnetic Hamiltonian system.\\

In this section, we first derive precisely the geometric constraint conditions of
the magnetic symplectic form for the dynamical vector field
of a magnetic Hamiltonian system, that is,
Type I and Type II of Hamilton-Jacobi equation for
the magnetic Hamiltonian system.
In order to do this, in the following we first give
some important notions and prove a key lemma, which is an important
tool for the proofs of two types of
Hamilton-Jacobi theorem for the magnetic Hamiltonian system.\\

Denote by $\Omega^i(Q)$ the set of all i-forms on $Q$, $i=1,2.$
For any $\gamma \in \Omega^1(Q),\; q\in Q, $ then $\gamma(q)\in T_q^*Q, $
and we can define a map $\gamma: Q \rightarrow T^*Q, \; q \rightarrow (q, \gamma(q)).$
Hence we say often that the map $\gamma: Q
\rightarrow T^*Q$ is an one-form on $Q$. If the one-form $\gamma$ is closed,
then $\mathbf{d}\gamma(x,y)=0, \; \forall\; x, y \in TQ$.
Note that for any $v, w \in TT^* Q, $ we have that
$\mathbf{d}\gamma(T\pi_{Q}(v),T\pi_{Q}(w))=\pi^*(\mathbf{d}\gamma )(v, w)$
is a two-form on the cotangent bundle $T^*Q$, where
$\pi^*: T^*Q \rightarrow T^*T^*Q.$ Thus,
in the following we can give a weaker notion.
\begin{defi}
The one-form $\gamma$ is called to be closed with respect to $T\pi_{Q}:
TT^* Q \rightarrow TQ, $ if for any $v, w \in TT^* Q, $ we have
that $\mathbf{d}\gamma(T\pi_{Q}(v),T\pi_{Q}(w))=0. $
\end{defi}

For the one-form $\gamma: Q \rightarrow T^*Q$, $\mathbf{d}\gamma$
is a two-form on $Q$. Assume that $B$ is a closed two-form on $Q$,
we say that the $\gamma$ satisfies condition $\mathbf{d}\gamma=-B$,
if for any $ x, y \in TQ$, we have that $(\mathbf{d}\gamma +B)(x,y)=0.$
In the following we can give a new notion.
\begin{defi}
Assume that $\gamma: Q
\rightarrow T^*Q$ is an one-form on $Q$,
we say that the $\gamma$ satisfies condition that
$\mathbf{d}\gamma=-B$ with respect to $T\pi_{Q}:
TT^* Q \rightarrow TQ, $ if for any $v, w \in TT^* Q, $ we have
that $(\mathbf{d}\gamma +B)(T\pi_{Q}(v),T\pi_{Q}(w))=0. $
\end{defi}

From the above definition we know that, if $\gamma$ satisfies condition
$\mathbf{d}\gamma=-B$, then it must satisfy condition
$\mathbf{d}\gamma=-B$ with respect to $T\pi_{Q}: TT^* Q \rightarrow
TQ. $ Conversely, if $\gamma$ satisfies condition
$\mathbf{d}\gamma=-B$ with respect to
$T\pi_{Q}: TT^* Q \rightarrow TQ, $ then it may not satisfy condition
$\mathbf{d}\gamma=-B$. We can
prove a general result as follows, which states that
the notion that $\gamma$ satisfies condition
$\mathbf{d}\gamma=-B$
with respect to $T\pi_{Q}: TT^* Q \rightarrow TQ, $
is not equivalent to the notion that $\gamma$ satisfies condition
$\mathbf{d}\gamma=-B$.

\begin{prop}
Assume that $\gamma: Q \rightarrow T^*Q$ is an one-form on $Q$ and
it doesn't satisfy condition
$\mathbf{d}\gamma=-B$. we define the set $N$, which is a subset of $TQ$,
such that the one-form $\gamma$ on $N$ satisfies the condition that
for any $x,y \in N, \; (\mathbf{d}\gamma +B)(x,y)\neq 0. $ Denote by
$Ker(T\pi_Q)= \{u \in TT^*Q| \; T\pi_Q(u)=0 \}, $ and $T\gamma: TQ
\rightarrow TT^* Q $ is the tangent map of $\gamma: Q \rightarrow T^*Q. $
If $T\gamma(N)\subset Ker(T\pi_Q), $ then
$\gamma$ satisfies condition
$\mathbf{d}\gamma=-B$ with respect to $T\pi_{Q}: TT^* Q \rightarrow TQ.
$\end{prop}

\noindent{\bf Proof: } In fact, for any $v, w \in TT^* Q, $ if
$T\pi_{Q}(v) \notin N, $ or $T\pi_{Q}(w))\notin N, $ then by the
definition of $N$, we know that
$(\mathbf{d}\gamma+B)(T\pi_{Q}(v),T\pi_{Q}(w))=0; $ If $T\pi_{Q}(v)\in
N, $ and $T\pi_{Q}(w))\in N, $ from the condition $T\gamma(N)\subset
Ker(T\pi_Q), $ we know that $T\pi_{Q}\cdot T\gamma \cdot
T\pi_{Q}(v)= T\pi_{Q}(v)=0, $ and $T\pi_{Q}\cdot T\gamma \cdot
T\pi_{Q}(w)= T\pi_{Q}(w)=0, $ where we have used the
relation $\pi_Q\cdot \gamma\cdot \pi_Q= \pi_Q, $ and hence
$(\mathbf{d}\gamma+B)(T\pi_{Q}(v),T\pi_{Q}(w))=0. $ Thus, for any $v, w
\in TT^* Q, $ we have always that
$(\mathbf{d}\gamma+B)(T\pi_{Q}(v),T\pi_{Q}(w))=0, $ that is, $\gamma$
satisfies condition $\mathbf{d}\gamma=-B$
with respect to $T\pi_{Q}: TT^* Q \rightarrow TQ. $
\hskip 0.3cm $\blacksquare$\\

From the above Definition 3.1 and Definition 3.2, we know that,
when $B=0$, the notion that, $\gamma$ satisfies condition
$\mathbf{d}\gamma=-B$ with respect to $T\pi_{Q}:
TT^* Q \rightarrow TQ, $ become the notion that
$\gamma$ is closed with respect to $T\pi_{Q}:
TT^* Q \rightarrow TQ. $  Now, we can prove the following lemma,
which  is a generalization of a corresponding to lemma given by
Wang \cite{wa17},
and the lemma is a very important tool for our research.

\begin{lemm}
Assume that $\gamma: Q \rightarrow T^*Q$ is an one-form on $Q$, and
$\lambda=\gamma \cdot \pi_{Q}: T^* Q \rightarrow T^* Q .$
For the magnetic symplectic form $\omega^B= \omega- \pi_Q^*B $ on $T^*Q$,
where $\omega$ is the canonical symplectic form on $T^*Q$,
and $B$ is a closed two-form on $Q$,
then we have that the following two assertions hold.\\
\noindent $(\mathrm{i})$ For any $v, w \in
TT^* Q, \; \lambda^*\omega^B(v,w)= -(\mathbf{d}\gamma+B)(T\pi_{Q}(v), \;
T\pi_{Q}(w))$; \\
\noindent $(\mathrm{ii})$ For any $v, w \in TT^* Q, \;
\omega^B(T\lambda \cdot v,w)= \omega^B(v, w-T\lambda \cdot
w)-(\mathbf{d}\gamma+B)(T\pi_{Q}(v), \; T\pi_{Q}(w)). $
\end{lemm}

\noindent{\bf Proof:} We first prove the assertion $(\mathrm{i})$.
Since $\omega$ is the canonical symplectic form on $T^*Q$,
we know that there is an unique canonical one-form $\theta$, such that
$\omega= -\mathbf{d} \theta. $ From the Proposition 3.2.11 in
Abraham and Marsden \cite{abma78}, we have that for the one-form
$\gamma: Q \rightarrow T^*Q, \; \gamma^* \theta= \gamma. $ Then we
can obtain that for any $x, y \in TQ,$
\begin{align*}
\gamma^*\omega(x,y) = \gamma^* (-\mathbf{d} \theta) (x, y) =
-\mathbf{d}(\gamma^* \theta)(x, y)= -\mathbf{d}\gamma (x, y).
\end{align*}
Note that $\lambda=\gamma \cdot \pi_{Q}: T^* Q \rightarrow T^* Q, $
and $\lambda^*= \pi_{Q}^* \cdot \gamma^*: T^*T^* Q \rightarrow
T^*T^* Q, $ then we have that  for any $v, w \in TT^* Q $,
\begin{align*}
\lambda^*\omega(v,w) &= \lambda^* (-\mathbf{d} \theta) (v, w)
=-\mathbf{d}(\lambda^* \theta)(v, w)= -\mathbf{d}(\pi_{Q}^* \cdot
\gamma^* \theta)(v, w)\\ &= -\mathbf{d}(\pi_{Q}^* \cdot\gamma )(v,
w)= -\mathbf{d}\gamma(T\pi_{Q}(v), \; T\pi_{Q}(w)).
\end{align*}
Hence, we have that
\begin{align*}
\lambda^*\omega^B(v,w)& =\lambda^*\omega(v,w)-\lambda^*\cdot \pi_Q^*B(v,w)\\
& =-\mathbf{d}\gamma(T\pi_{Q}(v), \; T\pi_{Q}(w))
-(\pi_Q\cdot \gamma \cdot \pi_{Q})^*B(v,w)\\
& =-\mathbf{d}\gamma(T\pi_{Q}(v), \; T\pi_{Q}(w))- \pi_Q^*B(v,w)\\
& =-(\mathbf{d}\gamma+B)(T\pi_{Q}(v), \; T\pi_{Q}(w)),
\end{align*}
where we have used the relation $\pi_Q\cdot \gamma\cdot \pi_Q= \pi_Q. $
It follows that the assertion $(\mathrm{i})$ holds.\\

Next, we prove the assertion $(\mathrm{ii})$. For any $v, w \in TT^*
Q,$ note that $v- T(\gamma \cdot \pi_Q)\cdot v$ is vertical, because
$$
T\pi_Q(v- T(\gamma \cdot \pi_Q)\cdot v)=T\pi_Q(v)-T(\pi_Q\cdot
\gamma\cdot \pi_Q)\cdot v= T\pi_Q(v)-T\pi_Q(v)=0,
$$
Thus, $\omega(v- T(\gamma \cdot \pi_Q)\cdot v,w- T(\gamma \cdot
\pi_Q)\cdot w)= 0, $ and hence,
$$\omega(T(\gamma \cdot \pi_Q)\cdot v, \; w)=
\omega(v, \; w-T(\gamma \cdot \pi_Q)\cdot w)+ \omega(T(\gamma \cdot
\pi_Q)\cdot v, \; T(\gamma \cdot \pi_Q)\cdot w). $$ However, the
second term on the right-hand side is given by
$$
\omega(T(\gamma \cdot \pi_Q)\cdot v, \; T(\gamma \cdot \pi_Q)\cdot
w)= \gamma^*\omega(T\pi_Q(v), \; T\pi_Q(w))=
-\mathbf{d}\gamma(T\pi_{Q}(v), \; T\pi_{Q}(w)),
$$
It follows that
\begin{align*}
\omega(T\lambda \cdot v,w) &=\omega(T(\gamma \cdot \pi_Q)\cdot v, \;
w)\\ &= \omega(v, \; w-T(\gamma \cdot \pi_Q)\cdot w)-\mathbf{d}\gamma(T\pi_{Q}(v), \; T\pi_{Q}(w))
\\ &= \omega(v,
w-T\lambda \cdot w)-\mathbf{d}\gamma(T\pi_{Q}(v), \; T\pi_{Q}(w)).
\end{align*}
Hence,  we have that
\begin{align*}
& \omega^B(T\lambda \cdot v,w)= \omega(T\lambda \cdot v,w)-\pi_Q^*B(T\lambda \cdot v,w)\\
& =\omega(v, w-T\lambda \cdot w)-\mathbf{d}\gamma(T\pi_{Q}(v), \; T\pi_{Q}(w))
-B(T\pi_Q\cdot T\lambda \cdot v, \; T\pi_{Q}(w))\\
& =\omega^B(v, w-T\lambda \cdot w)+\pi_Q^*B(v, w-T\lambda \cdot w)\\
& \;\;\;\;\; -\mathbf{d}\gamma(T\pi_{Q}(v), \; T\pi_{Q}(w))
-B(T(\pi_Q\cdot \lambda) \cdot v, \; T\pi_{Q}(w))\\
& =\omega^B(v, w-T\lambda \cdot w)+\pi_Q^*B(v, w)-B(T\pi_{Q}(v), \; T\pi_Q\cdot T\lambda \cdot w)\\
& \;\;\;\;\; -\mathbf{d}\gamma(T\pi_{Q}(v), \; T\pi_{Q}(w))
-B(T(\pi_Q\cdot \gamma \cdot \pi_{Q}) \cdot v, \; T\pi_{Q}(w))\\
& =\omega^B(v, w-T\lambda \cdot w)+\pi_Q^*B(v, w)-B(T\pi_{Q}(v), \; T(\pi_Q\cdot \lambda) \cdot w)\\
& \;\;\;\;\; -\mathbf{d}\gamma(T\pi_{Q}(v), \; T\pi_{Q}(w))
-B(T\pi_Q (v), \; T\pi_{Q}(w))\\
& =\omega^B(v, w-T\lambda \cdot w)+\pi_Q^*B(v, w)-B(T\pi_{Q}(v), \; T\pi_{Q}(w))
-(\mathbf{d}\gamma+B)(T\pi_{Q}(v), \; T\pi_{Q}(w))\\
& =\omega^B(v, w-T\lambda \cdot
w)-(\mathbf{d}\gamma+B)(T\pi_{Q}(v), \; T\pi_{Q}(w)).
\end{align*}
Thus, the assertion $(\mathrm{ii})$ holds.
\hskip 0.3cm $\blacksquare$\\

Since a magnetic Hamiltonian system is a Hamiltonian system defined by the
magnetic symplectic form, and it is a canonical Hamiltonian system
coupling the action of a magnetic field $B$. Usually, under the impact of
magnetic term $\pi_Q^*B$, the magnetic symplectic form
$\omega^B=\omega- \pi_Q^*B $,
in general, is not the canonical symplectic form $\omega$ on $T^*Q$,
we cannot prove the Hamilton-Jacobi theorem for the magnetic Hamiltonian system
just like same as the above Theorem 1.1. But,
in the following we can give precisely the geometric constraint conditions of
magnetic symplectic form for the dynamical vector field of
the magnetic Hamiltonian system, that is, Type I and Type II of
Hamilton-Jacobi equation for the magnetic Hamiltonian system.
At first, for a given magnetic Hamiltonian system
$(T^*Q,\omega^B,H)$ on $T^*Q$, by using
the above Lemma 3.4, magnetic symplectic form $\omega^B$
and the magnetic Hamiltonian vector field $X^B_{H}$,
we can derive precisely the following type I of
Hamilton-Jacobi equation for the magnetic Hamiltonian system $(T^*Q,\omega^B,H)$.

\begin{theo} (Type I of Hamilton-Jacobi Theorem for a Magnetic Hamiltonian System)
For a given magnetic Hamiltonian system $(T^*Q,\omega^B,H)$ with
the magnetic symplectic form $\omega^B= \omega- \pi_Q^*B $ on $T^*Q$,
where $\omega$ is the canonical symplectic form on $T^* Q$
and $B$ is a closed two-form on $Q$,
assume that $\gamma: Q
\rightarrow T^*Q$ is an one-form on $Q$, and
$X^\gamma = T\pi_{Q}\cdot X^B_H \cdot \gamma$,
where $X^B_H$ is the dynamical vector field
of the magnetic Hamiltonian system $(T^*Q,\omega^B,H)$,
that is, the magnetic Hamiltonian vector field.
If the one-form $\gamma: Q \rightarrow T^*Q $ satisfies the condition that
$\mathbf{d}\gamma=-B $ with respect to $T\pi_{Q}:
TT^* Q \rightarrow TQ, $ then $\gamma$ is a solution of the equation
$T\gamma\cdot X^\gamma= X^B_H\cdot \gamma .$
The equation is called the Type I of
Hamilton-Jacobi equation for the magnetic Hamiltonian system
$(T^*Q,\omega^B,H)$.
Here the maps involved in the theorem are shown
in the following Diagram-1.
\begin{center}
\hskip 0cm \xymatrix{ & T^* Q \ar[r]^{\pi_Q}
 & Q \ar[d]_{X^\gamma} \ar[r]^{\gamma} & T^*Q \ar[d]^{X^B_H} \\
 & T(T^*Q) & TQ \ar[l]_{T\gamma} & T(T^* Q)\ar[l]_{T\pi_Q}}
\end{center}
$$\mbox{Diagram-1}$$
\end{theo}
\noindent{\bf Proof: } If we take that
$v= X^B_H\cdot \gamma \in TT^* Q, $ and for
any $w \in TT^* Q, \; T\pi_{Q}(w)\neq 0, $ from Lemma 3.4(ii) and
$\mathbf{d}\gamma=-B $ with respect to $T\pi_{Q}:
TT^* Q \rightarrow TQ, $ that is,
$(\mathbf{d}\gamma+B)(T\pi_{Q}\cdot X^B_{H}\cdot\gamma, \; T\pi_{Q}\cdot w)=0,$
we have that
\begin{align*}
\omega^B(T\gamma \cdot X^\gamma, \; w)&
=\omega^B(T\gamma \cdot T\pi_{Q} \cdot X^B_H\cdot\gamma, \; w)
= \omega^B(T(\gamma \cdot \pi_Q)\cdot X^B_H\cdot \gamma, \; w)\\
&= \omega^B(X^B_H\cdot \gamma, \; w-T(\gamma \cdot \pi_Q)\cdot w)
-(\mathbf{d}\gamma+B)(T\pi_{Q}\cdot X^B_{H}\cdot\gamma, \; T\pi_{Q}\cdot w)\\
& = \omega^B(X^B_H\cdot \gamma, \; w) - \omega^B(X^B_H\cdot \gamma, \;
T\lambda \cdot w).
\end{align*}
Hence, we have that
\begin{equation}
\omega^B(T\gamma \cdot X^\gamma, \; w)- \omega^B(X^B_H\cdot \gamma, \; w)
= -\omega^B(X^B_H\cdot \gamma, \; T\lambda \cdot w). \; \label{3.1}
\end{equation}
If $\gamma$ satisfies the equation $T\gamma\cdot X^\gamma= X^B_H\cdot \gamma ,$
from Lemma 3.4(i) we know that the right side of (3.1) becomes that
\begin{align*}
\omega^B(X^B_H\cdot \gamma, \; T\lambda \cdot w) &
= \omega^B(T\gamma \cdot X^\gamma, \; T\lambda \cdot w)\\
&= \omega^B(T\gamma \cdot T\pi_{Q} \cdot X^B_H\cdot\gamma, \; T\lambda \cdot w)\\
&= \omega^B(T\lambda \cdot X^B_{H}\cdot\gamma, \; T\lambda \cdot w)\\
&= \lambda^*\omega^B(X^B_{H}\cdot\gamma, \; w)\\
&=-(\mathbf{d}\gamma+B)(T\pi_{Q}\cdot X^B_{H}\cdot\gamma, \; T\pi_{Q}\cdot w)=0,
\end{align*}
since $\gamma: Q \rightarrow T^*Q $ satisfies the condition that
$\mathbf{d}\gamma=-B $ with respect to $T\pi_{Q}:
TT^* Q \rightarrow TQ. $
But, because the magnetic symplectic form $\omega^B$ is non-degenerate,
the left side of (3.1) equals zero, only when
$\gamma$ satisfies the equation $T\gamma \cdot X^\gamma= X^B_H\cdot \gamma .$ Thus,
if the one-form $\gamma: Q \rightarrow T^*Q $ satisfies the condition that
$\mathbf{d}\gamma=-B$ with respect to $T\pi_{Q}:
TT^* Q \rightarrow TQ, $ then $\gamma$ must be a solution of
the Type I of Hamilton-Jacobi equation
$T\gamma\cdot X^\gamma= X^B_H\cdot \gamma ,$ for
the magnetic Hamiltonian system $(T^*Q,\omega^B,H)$.
\hskip 0.3cm $\blacksquare$\\

It is worthy of noting that, when $B=0$, in this case the magnetic symplectic form $\omega^B$
is just the canonical symplectic form $\omega$ on $T^*Q$,
and the magnetic Hamiltonian system $(T^*Q,\omega^B,H)$ becomes
the Hamiltonian system $(T^*Q,\omega,H)$ with the canonical symplectic form $\omega$,
and the condition that the one-form $\gamma: Q \rightarrow T^*Q $ satisfies the condition,
$\mathbf{d}\gamma=-B $ with respect to $T\pi_{Q}:
TT^* Q \rightarrow TQ, $ becomes the condition that  $\gamma $ is closed with respect to
$T\pi_Q: TT^* Q \rightarrow TQ.$ Thus, from above Theorem 3.5,
we can obtain Theorem 2.5 in Wang \cite{wa17}, that is.
the Type I of Hamilton-Jacobi theorem for a Hamiltonian system.
On the other hand, from the proof of Theorem 2.5 in Wang \cite{wa17},
we know that if an one-form $\gamma: Q \rightarrow T^*Q $ is not
closed on $Q$ with respect to $T\pi_Q: TT^* Q \rightarrow TQ $,
then $\gamma$ is not a solution of
the Type I of Hamilton-Jacobi equation
$T\gamma\cdot X^\gamma= X_H\cdot \gamma .$ But, note that,
 if $\gamma: Q \rightarrow T^*Q $ is not closed on $Q$ with respect to
$T\pi_Q: TT^* Q \rightarrow TQ, $ that is, there exist $v,w \in TT^*Q,$
such that $\mathbf{d}\gamma(T\pi_Q(v), \; T\pi_Q(w))\neq 0,$
and hence $\gamma$ is not yet closed on $Q$. However,
because $\mathbf{d}\cdot \mathbf{d}\gamma= \mathbf{d}^2 \gamma =0, $
and hence the $\mathbf{d}\gamma$ is a closed two-form on $Q$.
Thus, we can construct a magnetic symplectic form on $T^*Q$,
that is, $\omega^B= \omega+ \pi_Q^*(\mathbf{d}\gamma)=\omega- \pi_Q^*B , $
where $B=- \mathbf{d}\gamma,$ and
$\omega$ is the canonical symplectic form on $T^*Q$,
and $\pi_Q^*: T^*Q \rightarrow T^*T^*Q $. Moreover,
we hope to look for a new magnetic Hamiltonian system,
such that $\gamma$ is a solution of the Type I
of Hamilton-Jacobi equation for the new magnetic Hamiltonian system.
In fact, for a given Hamiltonian system $(T^*Q, \omega, H )$ with
the canonical symplectic form $\omega$ on $T^*Q$, and
$\gamma: Q \rightarrow T^*Q $ is an
one-form on $Q$, and it is not closed with respect to
$T\pi_Q: TT^* Q \rightarrow TQ. $ Then we can construct
a magnetic symplectic form on $T^*Q$,
$\omega^B= \omega+ \pi_Q^*(\mathbf{d}\gamma), $ where $B=- \mathbf{d}\gamma,$
and a magnetic Hamiltonian system $(T^*Q, \omega^B, H )$,
its dynamical vector field is given by
$X^B_{H}$, which satisfies the magnetic Hamiltonian equation, that is,
$\mathbf{i}_{X^B_{H} }\omega^B= \mathbf{d}H $.
In this case, for any $x, y \in TQ, $ we have that $(\mathbf{d}\gamma +B)(x, y)=0$,
and hence for any $v, w \in TT^* Q, $ we have
$(\mathbf{d}\gamma +B)(T\pi_{Q}(v),T\pi_{Q}(w))=0, $
that is, the one-form $\gamma: Q \rightarrow T^*Q $ satisfies the condition,
$\mathbf{d}\gamma=-B$ with respect to $T\pi_{Q}:
TT^* Q \rightarrow TQ. $ Thus,
by using Lemma 3.4 and the magnetic Hamiltonian vector field
$X^B_H$, from Theorem 3.5
we can obtain the following Theorem 3.6.
\begin{theo}
For a given Hamiltonian system $(T^*Q,\omega, H)$ with
the canonical symplectic form $\omega$ on $T^*Q$,
and assume that the one-form $\gamma: Q
\rightarrow T^*Q$ is not closed with respect to
$T\pi_Q: TT^* Q \rightarrow TQ. $ Then one can construct
a magnetic symplectic form on $T^*Q$, that is,
$\omega^B= \omega+ \pi_Q^*(\mathbf{d}\gamma), $
where $B=- \mathbf{d}\gamma,$
and a magnetic Hamiltonian system $(T^*Q, \omega^B, H )$.
Denote $X^\gamma = T\pi_{Q}\cdot X^B_H \cdot \gamma$,
where $X^B_H$ is the dynamical vector field
of the magnetic Hamiltonian system $(T^*Q,\omega^B,H)$.
Then $\gamma$ is a solution of the Type I of
Hamilton-Jacobi equation $T\gamma\cdot X^\gamma= X^B_H\cdot \gamma ,$
for the magnetic Hamiltonian system $(T^*Q,\omega^B,H)$.
\end{theo}

Next, for any symplectic map $\varepsilon: T^* Q \rightarrow T^* Q $
with respect to the magnetic symplectic form $\omega^B$,
we can derive precisely the following Type II of
Hamilton-Jacobi equation for the magnetic Hamiltonian system
$(T^*Q,\omega^B,H)$. For convenience,
the maps involved in the following theorem and its proof are shown
in Diagram-2.
\begin{center}
\hskip 0cm \xymatrix{ & T^* Q \ar[r]^{\varepsilon}
& T^*Q \ar[d]_{X^B_{H\cdot \varepsilon}}
\ar[dr]^{X^\varepsilon} \ar[r]^{\pi_Q}
& Q \ar[r]^{\gamma} & T^*Q \ar[d]^{X^B_H} \\
&  & T(T^*Q) & TQ \ar[l]_{T\gamma} & T(T^* Q)\ar[l]_{T\pi_Q}}
\end{center}
$$\mbox{Diagram-2}$$

\begin{theo}
(Type II of Hamilton-Jacobi Theorem for a Magnetic Hamiltonian System)
For the magnetic Hamiltonian system $(T^*Q,\omega^B,H)$ with the
magnetic symplectic form $\omega^B= \omega- \pi_Q^*B $ on $T^*Q$,
where $\omega$ is the canonical symplectic form on $T^* Q$
and $B$ is a closed two-form on $Q$,
assume that $\gamma: Q \rightarrow T^*Q$ is an one-form on $Q$, and
$\lambda=\gamma \cdot \pi_{Q}: T^* Q \rightarrow T^* Q $, and for any
symplectic map $\varepsilon: T^* Q \rightarrow T^* Q $ with respect to $\omega^B$,
denote by $ X^\varepsilon = T\pi_{Q}\cdot X^B_H \cdot \varepsilon$,
where $X^B_H$ is the dynamical vector field of the magnetic Hamiltonian system
$(T^*Q,\omega^B,H)$, that is, the magnetic Hamiltonian vector field.
Then $\varepsilon$ is a solution of the equation
$T\varepsilon\cdot X^B_{H\cdot\varepsilon}= T\lambda \cdot X^B_H \cdot \varepsilon,$
if and only if it is a solution of the equation $T\gamma \cdot X^\varepsilon= X^B_H\cdot
\varepsilon, $ where $ X^B_{H\cdot\varepsilon} \in
TT^*Q $ is the magnetic Hamiltonian vector field of the function $H\cdot\varepsilon:
T^*Q\rightarrow \mathbb{R} $.
The equation $T\gamma\cdot X^\varepsilon= X^B_H\cdot
\varepsilon ,$ is called the Type II of Hamilton-Jacobi equation
for the magnetic Hamiltonian system $(T^*Q,\omega^B,H)$.
\end{theo}
\noindent{\bf Proof: }
If we take that $v= X^B_H\cdot \varepsilon \in TT^* Q, $ and for
any $w \in TT^* Q, \; T\lambda(w)\neq 0, $ from Lemma 3.4 we have that
\begin{align*}
&\omega^B(T\gamma \cdot X^\varepsilon, \; w)
= \omega^B(T\gamma \cdot T\pi_Q\cdot X^B_H\cdot \varepsilon, \; w)
= \omega^B(T(\gamma \cdot \pi_Q)\cdot X^B_H\cdot \varepsilon, \; w)\\
&= \omega^B(X^B_H\cdot \varepsilon, \; w-T(\gamma \cdot \pi_Q)\cdot w)
-(\mathbf{d}\gamma+B)(T\pi_{Q}(X^B_H\cdot \varepsilon), \; T\pi_{Q}(w))\\
& =\omega^B(X^B_H\cdot \varepsilon, \; w) - \omega^B(X^B_H\cdot \varepsilon, \;
T\lambda \cdot w)+\lambda^*\omega^B(X^B_H\cdot \varepsilon, \; w)\\
& =\omega^B(X^B_H\cdot \varepsilon, \; w) - \omega^B(X^B_H\cdot \varepsilon, \;
T\lambda \cdot w)+ \omega^B(T\lambda \cdot X^B_H\cdot \varepsilon, \; T\lambda \cdot w).
\end{align*}
Because $\varepsilon: T^* Q
\rightarrow T^* Q $ is symplectic with respect to $\omega^B$,
and hence $ X^B_H\cdot \varepsilon= T\varepsilon \cdot X^B_{H\cdot\varepsilon}, $
along $\varepsilon$.
From the above arguments, we can obtain that
\begin{align*}
&\omega^B(T\gamma \cdot X^\varepsilon, \; w)- \omega^B(X^B_H\cdot \varepsilon, \; w)\\
& =- \omega^B(X^B_H\cdot \varepsilon, \; T\lambda \cdot w)
+ \omega^B(T\lambda \cdot X^B_H\cdot \varepsilon, \; T\lambda \cdot w)\\
& =-\omega^B(T\varepsilon \cdot X^B_{H\cdot\varepsilon}, \; T\lambda \cdot w)
+ \omega^B(T\lambda \cdot X^B_H \cdot \varepsilon, \; T\lambda \cdot w)\\
& = \omega^B(T\lambda \cdot X^B_H \cdot \varepsilon
-T\varepsilon \cdot X^B_{H\cdot\varepsilon}, \; T\lambda \cdot w).
\end{align*}
Because the magnetic symplectic form $\omega^B$ is non-degenerate,
it follows that $T\gamma \cdot X^\varepsilon= X^B_H\cdot
\varepsilon ,$ is equivalent to $T\varepsilon \cdot X^B_{H\cdot\varepsilon}
= T\lambda\cdot X^B_H\cdot \varepsilon $.
Thus, $\varepsilon$ is a solution of the equation
$T\varepsilon\cdot X^B_{H\cdot\varepsilon}= T\lambda \cdot X^B_H \cdot \varepsilon,$
if and only if it is a solution of the Type II of Hamilton-Jacobi equation
$T\gamma \cdot X^\varepsilon= X^B_H \cdot \varepsilon .$
\hskip 0.3cm $\blacksquare$

\begin{rema}
It is worthy of noting that, the Type I of Hamilton-Jacobi equation
$T\gamma \cdot X^\gamma= X^B_H \cdot \gamma ,$
is the equation of the differential one-form $\gamma$; and
the Type II of Hamilton-Jacobi equation $T\gamma\cdot X^\varepsilon
= X^B_H \cdot \varepsilon ,$ is the equation of
the symplectic diffeomorphism map $\varepsilon$.
When $B=0$, in this case the magnetic symplectic form $\omega^B$
is just the canonical symplectic form $\omega$ on $T^*Q$, and
the magnetic Hamiltonian system is just the canonical Hamiltonian system itself.
From the above Type I and Type II of Hamilton-Jacobi theorems, that is,
Theorem 3.5 and Theorem 3.7, we can get the Theorem 2.5
and Theorem 2.6 in Wang \cite{wa17}.
It shows that Theorem 3.5 and Theorem 3.7 can be regarded as an extension of two types of
Hamilton-Jacobi theorem for Hamiltonian system given in \cite{wa17} to that for the
magnetic Hamiltonian system.
\end{rema}

\section{Hamilton-Jacobi Equation for Distributional Magnetic Hamiltonian System }

In this section we shall derive precisely
the geometric constraint conditions of the induced distributional two-form
for the nonholonomic dynamical vector field of
distributional magnetic Hamiltonian system,
that is, the two types of Hamilton-Jacobi equation
for the distributional magnetic Hamiltonian system.
In order to do this, in the following we first give
some important notions and prove a key lemma, which is an important
tool for the proofs of two types of
Hamilton-Jacobi theorem for the distributional magnetic Hamiltonian system.\\

Assume that $\mathcal{D}\subset TQ$ is a $\mathcal{D}$-regularly nonholonomic
constraint, and the constraint submanifold
$\mathcal{M}=\mathcal{F}L(\mathcal{D})\subset T^*Q$,
the distribution
$\mathcal{F}=(T\pi_Q)^{-1}(\mathcal{D})\subset TT^*Q,$ and
$\gamma: Q \rightarrow T^*Q$ is an one-form on $Q$, and
$B$ is a closed two-form on $Q$, in the following we first introduce two weaker notions.

\begin{defi}
\noindent $(\mathrm{i})$ The one-form $\gamma$ is called to be closed
on $\mathcal{D}$ with respect to $T\pi_{Q}:
TT^* Q \rightarrow TQ, $ if for any $v, w \in TT^* Q, $
and $T\pi_{Q}(v), \; T\pi_{Q}(w) \in \mathcal{D},$  we have
that $\mathbf{d}\gamma(T\pi_{Q}(v),T\pi_{Q}(w))=0; $\\

\noindent $(\mathrm{ii})$ The one-form $\gamma: Q
\rightarrow T^*Q$ is called that satisfies condition that
$\mathbf{d}\gamma=-B$ on $\mathcal{D}$ with respect to $T\pi_{Q}:
TT^* Q \rightarrow TQ, $ if for any $v, w \in TT^* Q, $
and $T\pi_{Q}(v), \; T\pi_{Q}(w) \in \mathcal{D},$ we have
that  $(\mathbf{d}\gamma +B)(T\pi_{Q}(v),T\pi_{Q}(w))=0. $
\end{defi}

From the above Definition 4.1, we know that,
when $B=0$, the notion that, $\gamma$ satisfies condition
that $\mathbf{d}\gamma=-B$ on $\mathcal{D}$ with respect to $T\pi_{Q}:
TT^* Q \rightarrow TQ, $ become the notion that
$\gamma$ is closed on $\mathcal{D}$ with respect to $T\pi_{Q}:
TT^* Q \rightarrow TQ. $ On the other hand, it is worthy of noting that
the notion that $\gamma$ satisfies condition
that $\mathbf{d}\gamma=-B$ on $\mathcal{D}$ with respect to $T\pi_{Q}:
TT^* Q \rightarrow TQ, $ is weaker than the notion that $\gamma$
satisfies condition $\mathbf{d}\gamma=-B$ on $\mathcal{D},$
that is, $(\mathbf{d}\gamma+B)(x,y)=0, \; \forall\; x, y \in \mathcal{D}$.
In fact, if $\gamma$  satisfies condition
$\mathbf{d}\gamma=-B$ on $\mathcal{D}$,
then it must satisfy condition that $\mathbf{d}\gamma=-B$
on $\mathcal{D}$ with respect to $T\pi_{Q}: TT^* Q \rightarrow TQ. $
Conversely, if $\gamma$  satisfies condition
that $\mathbf{d}\gamma=-B$ on $\mathcal{D}$ with respect to
$T\pi_{Q}: TT^* Q \rightarrow TQ, $ then it may not satisfy condition
$\mathbf{d}\gamma=-B$ on $\mathcal{D}$.
We can prove a general result as follows, which states that
the notion that, the $\gamma$ satisfies condition that
$\mathbf{d}\gamma=-B$ on $\mathcal{D}$
with respect to $T\pi_{Q}: TT^* Q \rightarrow TQ, $
is not equivalent to the notion that $\gamma$ satisfies condition
$\mathbf{d}\gamma=-B$ on $\mathcal{D}$.

\begin{prop}
Assume that $\gamma: Q \rightarrow T^*Q$ is an one-form on $Q$ and
it doesn't satisfy condition $\mathbf{d}\gamma=-B$ on $\mathcal{D}$.
We define the set $N$, which is a subset of $TQ$,
such that the one-form $\gamma$ on $N$ satisfies the condition that
for any $x,y \in N, \; (\mathbf{d}\gamma+B)(x,y)\neq 0. $ Denote
$Ker(T\pi_Q)= \{u \in TT^*Q| \; T\pi_Q(u)=0 \}, $ and $T\gamma: TQ
\rightarrow TT^* Q .$ If $T\gamma(N)\subset Ker(T\pi_Q), $ then
$\gamma$ satisfies condition $\mathbf{d}\gamma=-B$
with respect to $T\pi_{Q}: TT^* Q \rightarrow TQ.$
and hence $\gamma$ satisfies condition $\mathbf{d}\gamma=-B$
on $\mathcal{D}$ with respect to
$T\pi_{Q}: TT^* Q \rightarrow TQ.$
\end{prop}

\noindent{\bf Proof: } If the $\gamma: Q \rightarrow T^*Q$
doesn't satisfy condition $\mathbf{d}\gamma=-B$ on $\mathcal{D}$,
then it doesn't yet satisfy condition $\mathbf{d}\gamma=-B$.
From the proof of Lemma 3.3, for any $v, w
\in TT^* Q, $ we have always that
$(\mathbf{d}\gamma+B)(T\pi_{Q}(v),T\pi_{Q}(w))=0. $
In particular, for any $v, w \in TT^* Q, $
and $T\pi_{Q}(v), \; T\pi_{Q}(w) \in \mathcal{D},$  we have
$(\mathbf{d}\gamma+B)(T\pi_{Q}(v),T\pi_{Q}(w))=0. $
that is, $\gamma$ satisfies condition that $\mathbf{d}\gamma=-B$ on $\mathcal{D}$
with respect to $T\pi_{Q}: TT^* Q \rightarrow TQ. $
\hskip 0.3cm $\blacksquare$\\

Now, we prove the following Lemma 4.3. It is worthy of noting that
this lemma and Lemma 3.4 given in \S3 are the
important tool for the proofs of the two types of Hamilton-Jacobi
theorems for the distributional magnetic Hamiltonian system and the nonholonomic
reduced distributional magnetic Hamiltonian system.

\begin{lemm}
Assume that $\gamma: Q \rightarrow T^*Q$ is an one-form on $Q$, and
$\lambda=\gamma \cdot \pi_{Q}: T^* Q \rightarrow T^* Q ,$ and
$\omega$ is the canonical symplectic form on $T^*Q$, and
$\omega^B= \omega- \pi_Q^*B $
is the magnetic symplectic form on $T^*Q$.
If the Lagrangian $L$ is $\mathcal{D}$-regular, and
$\textmd{Im}(\gamma)\subset \mathcal{M}=\mathcal{F}L(\mathcal{D}), $
then we have that $ X^B_{H}\cdot \gamma \in \mathcal{F}$ along
$\gamma$, and $ X^B_{H}\cdot \lambda \in \mathcal{F}$ along
$\lambda$, that is, $T\pi_{Q}(X^B_H\cdot\gamma(q))\in
\mathcal{D}_{q}, \; \forall q \in Q $, and $T\pi_{Q}(X^B_H\cdot\lambda(q,p))\in
\mathcal{D}_{q}, \; \forall q \in Q, \; (q,p) \in T^* Q. $
Moreover, if a symplectic map $\varepsilon: T^* Q \rightarrow T^* Q $
with respect to the magnetic symplectic form $\omega^B$ satisfies the
condition $\varepsilon(\mathcal{M})\subset \mathcal{M},$ then
we have that $ X^B_{H}\cdot \varepsilon \in \mathcal{F}$ along
$\varepsilon. $
\end{lemm}

\noindent{\bf Proof:} Under the canonical cotangent bundle coordinates, for any $q \in Q
, \; (q,p)\in T^* Q, $ we have that
$$
X^B_H\cdot \gamma(q)= (\sum^n_{i=1} (\frac{\partial H}{\partial
p_i}\frac{\partial}{\partial q^i} - \frac{\partial H}{\partial
q^i}\frac{\partial}{\partial p_i})-
\sum^n_{i,j=1}B_{ij}\frac{\partial H}{\partial
p_j}\frac{\partial}{\partial p_i})\gamma(q).
$$
and
$$
X^B_H\cdot \lambda(q,p)= (\sum^n_{i=1} (\frac{\partial H}{\partial
p_i}\frac{\partial}{\partial q^i} - \frac{\partial H}{\partial
q^i}\frac{\partial}{\partial p_i})-
\sum^n_{i,j=1}B_{ij}\frac{\partial H}{\partial
p_j}\frac{\partial}{\partial p_i})\gamma\cdot \pi_Q(q,p).
$$
Then,
$$
T\pi_Q(X^B_H\cdot \gamma(q))=T\pi_Q(X^B_H\cdot \lambda(q,p))
=\sum^n_{i=1}(\frac{\partial H}{\partial
p_i}\frac{\partial}{\partial q^i})\gamma(q) \in T_q Q.
$$
Since $\textmd{Im}(\gamma)\subset \mathcal{M}, $ and
$\gamma(q)\in \mathcal{M}_{(q,p)}=\mathcal{F}L(\mathcal{D}_q), $ from the Lagrangian $L$ is
$\mathcal{D}$-regular, and $\mathcal{F}L$ is a diffeomorphism, then
there exists a point $(q,\; v_q)\in \mathcal{D}_q, $ such that
$\mathcal{F}L(q,\; v_q)=\gamma(q). $ Thus,
$$
T\pi_Q(X^B_H\cdot \gamma(q))=T\pi_Q(X^B_H\cdot
\lambda(q,p))=\mathcal{F}L(q, \; v_q)\sum^n_{i=1}(\frac{\partial H}{\partial
p_i}\frac{\partial}{\partial q^i}) \in
\mathcal{D}_q,
$$
it follows that $ X^B_{H}\cdot \gamma \in \mathcal{F}$ along
$\gamma$, and $ X^B_{H}\cdot \lambda \in \mathcal{F}$ along
$\lambda$. Moreover, for the symplectic map $\varepsilon: T^* Q \rightarrow T^* Q $
with respect to the magnetic symplectic form $\omega^B$, we have that
$$
X^B_H\cdot \varepsilon (q,p)= (\sum^n_{i=1} (\frac{\partial H}{\partial
p_i}\frac{\partial}{\partial q^i} - \frac{\partial H}{\partial
q^i}\frac{\partial}{\partial p_i})-
\sum^n_{i,j=1}B_{ij}\frac{\partial H}{\partial
p_j}\frac{\partial}{\partial p_i})\varepsilon (q,p).
$$
If $\varepsilon$ satisfies the
condition $\varepsilon(\mathcal{M})\subset \mathcal{M},$
then for any $(q,p)\in \mathcal{M}_{(q,p)}$, we have that
$\varepsilon(q,p)\in \mathcal{M}_{(q,p)},$
and there exists a point $(q,\; v_q)\in \mathcal{D}_q, $ such that
$\mathcal{F}L(q,\; v_q)=\varepsilon (q,p). $ Thus,
$$
T\pi_Q(X^B_H\cdot \varepsilon(q,p))
= \sum^n_{i=1}(\frac{\partial H}{\partial p_i}\frac{\partial}{\partial q^i})
\varepsilon(q, p)
=\mathcal{F}L(q, \; v_q)\sum^n_{i=1}(\frac{\partial H}{\partial
p_i}\frac{\partial}{\partial q^i}) \in
\mathcal{D}_q,
$$
it follows that $ X^B_{H}\cdot \varepsilon \in \mathcal{F}$ along
$\varepsilon$.
\hskip 0.3cm $\blacksquare$\\

We note that for a nonholonomic magnetic Hamiltonian system,
under the restriction given by nonholonomic constraint,
in general, the dynamical vector field of a nonholonomic magnetic
Hamiltonian system may not be Hamiltonian.
On the other hand, since the distributional magnetic
Hamiltonian system is determined by a non-degenerate distributional two-form
induced from the magnetic symplectic form, but, the non-degenerate distributional two-form
is not a "true two-form" on a manifold, and hence the leading
distributional magnetic Hamiltonian system can not be Hamiltonian.
Thus, we can not describe the Hamilton-Jacobi equations for the nonholonomic
magnetic Hamiltonian system from the viewpoint of generating function
as in the classical Hamiltonian case, that is,
we cannot prove the Hamilton-Jacobi
theorem for the nonholonomic magnetic Hamiltonian system,
just like same as the above Theorem 1.1.
Since the distributional magnetic Hamiltonian system is a
dynamical system closely related to a magnetic Hamiltonian system,
in the following by using Lemma 3.4, Lemma 4.3,
and the non-degenerate
distributional two-form $\omega^B_{\mathcal{K}}$ and the
nonholonomic dynamical vector field $X^B_\mathcal {K}$ given
in \S 2 for the distributional magnetic Hamiltonian system,
we can derive precisely the geometric constraint conditions of
the non-degenerate distributional two-form $\omega^B_{\mathcal{K}}$
for the nonholonomic dynamical vector field $X^B_\mathcal {K}$,
that is, the two types of Hamilton-Jacobi equation for the distributional
magnetic Hamiltonian system $(\mathcal{K},\omega^B_{\mathcal {K}},H_{\mathcal{K}})$.
At first, we  prove the following Type I of
Hamilton-Jacobi theorem for the distributional magnetic Hamiltonian system.

\begin{theo} (Type I of Hamilton-Jacobi Theorem for the Distributional Magnetic Hamiltonian System)
For the nonholonomic magnetic Hamiltonian system $(T^*Q,\omega^B,\mathcal{D},H)$
with an associated distributional magnetic Hamiltonian system
$(\mathcal{K},\omega^B_{\mathcal {K}},H_{\mathcal{K}})$, assume that $\gamma: Q
\rightarrow T^*Q$ is an one-form on $Q$, and $X^\gamma =
T\pi_{Q}\cdot X^B_H \cdot \gamma$, where $X^B_{H}$ is
the magnetic Hamiltonian vector field of the associated
unconstrained magnetic Hamiltonian system
$(T^*Q,\omega^B,H)$.
Moreover, assume that $\textmd{Im}(\gamma)\subset
\mathcal{M}=\mathcal{F}L(\mathcal{D}), $ and $
\textmd{Im}(T\gamma)\subset \mathcal{K}. $ If the
one-form $\gamma: Q \rightarrow T^*Q $ satisfies the condition,
$\mathbf{d}\gamma=-B $ on $\mathcal{D}$ with respect to
$T\pi_Q: TT^* Q \rightarrow TQ, $ then $\gamma$ is a
solution of the equation $T\gamma \cdot
X^\gamma= X^B_{\mathcal{K}} \cdot \gamma. $ Here
$X^B_{\mathcal{K}}$ is the nonholonomic dynamical vector field
of the distributional magnetic Hamiltonian system
$(\mathcal{K},\omega^B_{\mathcal {K}},H_{\mathcal{K}})$. The equation $T\gamma \cdot
X^\gamma= X^B_{\mathcal{K}} \cdot \gamma $ is called the Type I of
Hamilton-Jacobi equation for the distributional magnetic Hamiltonian system
$(\mathcal{K},\omega^B_{\mathcal {K}},H_{\mathcal{K}})$. Here the
maps involved in the theorem are shown in
the following Diagram-3.
\begin{center}
\hskip 0cm \xymatrix{& \mathcal{M} \ar[d]_{X^B_{\mathcal{K}}}
\ar[r]^{i_{\mathcal{M}}} & T^* Q \ar[d]_{X^B_{H}}
 \ar[r]^{\pi_Q}
& Q \ar[d]_{{X}^{\gamma}} \ar[r]^{\gamma} & T^*Q \ar[d]^{X^B_H} \\
& \mathcal{K}  & T(T^*Q) \ar[l]_{\tau_{\mathcal{K}}} & TQ
\ar[l]_{T\gamma} & T(T^* Q)\ar[l]_{T\pi_Q}}
\end{center}
$$\mbox{Diagram-3}$$
\end{theo}

\noindent{\bf Proof: } At first, we note that
$\textmd{Im}(\gamma)\subset \mathcal{M}, $ and
$\textmd{Im}(T\gamma)\subset \mathcal{K}, $ in this case,
$\omega^B_{\mathcal{K}}\cdot
\tau_{\mathcal{K}}=\tau_{\mathcal{K}}\cdot \omega^B_{\mathcal{M}}=
\tau_{\mathcal{K}}\cdot i_{\mathcal{M}}^* \cdot \omega^B, $ along
$\textmd{Im}(T\gamma)$. Moreover, from the distributional magnetic Hamiltonian equation (2.2),
we have that $X^B_{\mathcal{K}}= \tau_{\mathcal{K}}\cdot X^B_H,$
and $\tau_{\mathcal{K}}\cdot X^B_{H}\cdot \gamma = X^B_{\mathcal{K}}\cdot \gamma $.
Thus, using the non-degenerate
distributional two-form $\omega^B_{\mathcal{K}}$, from Lemma 3.4(ii) and Lemma 4.3,
if we take that $v= X^B_{H}\cdot \gamma \in \mathcal{F},$ and for any $w
\in \mathcal{F}, \; T\lambda(w)\neq 0, $ and
$\tau_{\mathcal{K}}\cdot w \neq 0, $ then we have that
\begin{align*}
& \omega^B_{\mathcal{K}}(T\gamma \cdot X^\gamma, \;
\tau_{\mathcal{K}}\cdot w)=
\omega^B_{\mathcal{K}}(\tau_{\mathcal{K}}\cdot T\gamma \cdot
X^\gamma, \; \tau_{\mathcal{K}}\cdot w)\\ & =
\tau_{\mathcal{K}}\cdot i_{\mathcal{M}}^* \cdot \omega^B(T\gamma \cdot
T\pi_Q \cdot X^B_H \cdot \gamma, \; w ) = \tau_{\mathcal{K}}\cdot
i_{\mathcal{M}}^* \cdot \omega^B (T(\gamma \cdot \pi_Q)\cdot X^B_H \cdot \gamma, \; w)\\
& =\tau_{\mathcal{K}}\cdot i_{\mathcal{M}}^* \cdot
(\omega^B (X^B_H \cdot \gamma, \; w-T(\gamma \cdot \pi_Q)\cdot w)
-(\mathbf{d}\gamma+B)(T\pi_{Q}\cdot X^B_{H}\cdot\gamma, \; T\pi_{Q}\cdot w))\\
& = \tau_{\mathcal{K}}\cdot i_{\mathcal{M}}^* \cdot \omega^B (X^B_H \cdot
\gamma, \; w) - \tau_{\mathcal{K}}\cdot i_{\mathcal{M}}^* \cdot
\omega^B (X^B_H \cdot \gamma, \; T(\gamma
\cdot \pi_Q) \cdot w) \\
& \;\;\;\;\;\; - \tau_{\mathcal{K}}\cdot i_{\mathcal{M}}^* \cdot
(\mathbf{d}\gamma+B)(T\pi_{Q}\cdot X^B_{H}\cdot\gamma, \; T\pi_{Q}\cdot w)\\
& = \omega^B_{\mathcal{K}}( \tau_{\mathcal{K}}\cdot X^B_H \cdot \gamma,
\; \tau_{\mathcal{K}}\cdot w) -
\omega^B_{\mathcal{K}}(\tau_{\mathcal{K}}\cdot X^B_H \cdot \gamma, \;
\tau_{\mathcal{K}}\cdot T(\gamma \cdot \pi_Q) \cdot w)\\
& \;\;\;\;\;\; - \tau_{\mathcal{K}}\cdot i_{\mathcal{M}}^* \cdot
(\mathbf{d}\gamma+B)(T\pi_{Q}\cdot X^B_{H}\cdot\gamma, \; T\pi_{Q}\cdot w)\\
& = \omega^B_{\mathcal{K}}(X^B_{\mathcal{K}}\cdot \gamma, \;
\tau_{\mathcal{K}} \cdot w) -
\omega^B_{\mathcal{K}}(X^B_{\mathcal{K}} \cdot \gamma, \;
\tau_{\mathcal{K}}\cdot T\gamma \cdot T\pi_{Q}(w))\\
& \;\;\;\;\;\; - \tau_{\mathcal{K}}\cdot i_{\mathcal{M}}^* \cdot
(\mathbf{d}\gamma+B)(T\pi_{Q}\cdot X^B_{H}\cdot\gamma, \; T\pi_{Q}\cdot w),
\end{align*}
where we have used that $ \tau_{\mathcal{K}}\cdot T\gamma= T\gamma, $
since $\textmd{Im}(T\gamma)\subset \mathcal{K}, $ and
$\tau_{\mathcal{K}}\cdot X^B_H\cdot \gamma
= X^B_{\mathcal{K}}\cdot \gamma \in \mathcal{K}. $
Note that $X^B_{H}\cdot \gamma, \; w \in \mathcal{F},$ and
$T\pi_{Q}(X^B_H\cdot \gamma), \; T\pi_{Q}(w) \in \mathcal{D}. $
If the one-form $\gamma: Q \rightarrow T^*Q $ satisfies the condition,
$\mathbf{d}\gamma=-B $ on $\mathcal{D}$ with respect to
$T\pi_Q: TT^* Q \rightarrow TQ, $ then
$(\mathbf{d}\gamma+B)(T\pi_{Q}\cdot X^B_{H}\cdot\gamma, \; T\pi_{Q}\cdot w)=0,$
and hence
$$
\tau_{\mathcal{K}}\cdot i_{\mathcal{M}}^* \cdot(\mathbf{d}\gamma+B)
(T\pi_{Q}(X^B_H\cdot \gamma), \; T\pi_{Q}(w))=0,
$$
Thus, we have that
\begin{equation}
\omega^B_{\mathcal{K}}(T\gamma \cdot X^\gamma, \;
\tau_{\mathcal{K}}\cdot w)- \omega^B_{\mathcal{K}}(X^B_{\mathcal{K}}\cdot \gamma, \;
\tau_{\mathcal{K}} \cdot w)
= -\omega^B_{\mathcal{K}}(X^B_{\mathcal{K}}\cdot
\gamma, \; \tau_{\mathcal{K}}\cdot T\gamma \cdot T\pi_{Q}(w)).
\label{4.1} \end{equation}
If $\gamma$ satisfies the equation $T\gamma \cdot X^\gamma= X^B_{\mathcal{K}}\cdot \gamma ,$
from Lemma 3.4(i) we know that the right side of (4.1) becomes that
\begin{align*}
 -\omega^B_{\mathcal{K}}(X^B_{\mathcal{K}} \cdot \gamma, \;
\tau_{\mathcal{K}}\cdot T\gamma \cdot T\pi_{Q}(w))
& = -\omega^B_{\mathcal{K}}(T\gamma\cdot X^\gamma, \;
\tau_{\mathcal{K}}\cdot T\gamma \cdot T\pi_{Q}(w))\\
& = -\omega^B_{\mathcal{K}}(\tau_{\mathcal{K}}\cdot T\gamma\cdot X^\gamma, \;
\tau_{\mathcal{K}}\cdot T\gamma \cdot T\pi_{Q}(w))\\
& = -\tau_{\mathcal{K}}\cdot
i_{\mathcal{M}}^* \cdot \omega^B(T\gamma
\cdot T\pi_{Q}(X^B_{H}\cdot\gamma), \; T\gamma \cdot T\pi_{Q}(w))\\
& = -\tau_{\mathcal{K}}\cdot
i_{\mathcal{M}}^* \cdot \lambda^* \omega^B (X^B_{H}\cdot\gamma, \; w)\\
& = \tau_{\mathcal{K}}\cdot
i_{\mathcal{M}}^* \cdot (\mathbf{d}\gamma+B)(T\pi_{Q}\cdot X^B_{H}\cdot\gamma, \; T\pi_{Q}\cdot w)=0.
\end{align*}
Because the distributional two-form $\omega^B_{\mathcal{K}}$ is non-degenerate,
the left side of (4.1) equals zero, only when
$\gamma$ satisfies the equation $T\gamma\cdot X^\gamma= X^B_{\mathcal{K}}\cdot \gamma .$ Thus,
if the one-form $\gamma: Q \rightarrow T^*Q $ satisfies the condition that
$\mathbf{d}\gamma=-B $ on $\mathcal{D}$ with respect to
$T\pi_Q: TT^* Q \rightarrow TQ, $ then $\gamma$ must be a solution of the Type I of Hamilton-Jacobi equation
$T\gamma\cdot X^\gamma= X^B_{\mathcal{K}}\cdot \gamma ,$
for the distributional magnetic Hamiltonian system
$(\mathcal{K},\omega^B_{\mathcal {K}},H_{\mathcal{K}})$.
\hskip 0.3cm $\blacksquare$\\

It is worthy of noting that, when $B=0$, in this case the magnetic symplectic form $\omega^B$
is just the canonical symplectic form $\omega$ on $T^*Q$,
and the nonholonomic magnetic Hamiltonian system $(T^*Q,\omega^B,\mathcal{D},H)$ becomes
the nonholonomic Hamiltonian system $(T^*Q,\omega,\mathcal{D},H)$
with the canonical symplectic form $\omega$,
and the distributional magnetic Hamiltonian system
$(\mathcal{K},\omega^B_{\mathcal {K}},H_{\mathcal{K}})$ becomes the distributional Hamiltonian system
$(\mathcal{K},\omega_{\mathcal {K}},H_{\mathcal{K}})$,
and the condition that the one-form $\gamma: Q \rightarrow T^*Q $ satisfies the condition that
$\mathbf{d}\gamma=-B $ on $\mathcal{D}$ with respect to $T\pi_{Q}:
TT^* Q \rightarrow TQ, $ becomes that  $\gamma $ is closed on $\mathcal{D}$ with respect to
$T\pi_Q: TT^* Q \rightarrow TQ.$ Thus, from above Theorem 4.4,
we can obtain Theorem 3.5 in Le\'{o}n and Wang \cite{lewa15}, that is.
the Type I of Hamilton-Jacobi theorem for the distributional Hamiltonian system.
On the other hand, from the proofs of Theorem 3.5
in Le\'{o}n and Wang \cite{lewa15}, we know that, if the one-form
$\gamma: Q \rightarrow T^*Q $ is not closed on $\mathcal{D}$ with respect to
$T\pi_Q: TT^* Q \rightarrow TQ, $ then the $\gamma$ is not yet closed on $\mathcal{D}$,
that is,  $\mathbf{d}\gamma(x,y)\neq 0, \; \forall\;
x, y \in \mathcal{D}$, and hence $\gamma$ is not yet closed on $Q$.
However, in this case, we note that
$\mathbf{d}\cdot \mathbf{d}\gamma= \mathbf{d}^2 \gamma =0, $
and hence the $\mathbf{d}\gamma$ is a closed two-form on $Q$. Thus, we can construct
a magnetic symplectic form on $T^*Q$, $\omega^B= \omega+ \pi_Q^*(\mathbf{d}\gamma), $
where $B=- \mathbf{d}\gamma,$.
Moreover,  we can also construct a nonholonomic magnetic
Hamiltonian system $(T^*Q, \omega^B,\mathcal{D}, H )$
with an associated distributional magnetic Hamiltonian system
$(\mathcal{K},\omega^B_{\mathcal{K}},H_{\mathcal{K}})$, which satisfies
the distributional magnetic Hamiltonian equation (2.2),
$\mathbf{i}_{X^B_{\mathcal{K}} }\omega^B_{\mathcal{K}}= \mathbf{d}H_{\mathcal{K}}$.
In this case, the one-form $\gamma: Q \rightarrow T^*Q $ satisfies also the condition that
$\mathbf{d}\gamma=-B$ on $\mathcal{D}$ with respect to $T\pi_{Q}:
TT^* Q \rightarrow TQ, $
by using Lemma 3.4, Lemma 4.3, and the magnetic Hamiltonian vector field
$X^B_H$, from Theorem 4.4
we can obtain the following Theorem 4.5.
\begin{theo}
For a given nonholonomic Hamiltonian system $(T^*Q,\omega,\mathcal{D},H)$ with
the canonical symplectic form $\omega$ on $T^*Q$ and
$\mathcal{D}$-completely and $\mathcal{D}$-regularly nonholonomic
constraint $\mathcal{D} \subset TQ$,
and assume that the one-form $\gamma: Q
\rightarrow T^*Q$ is not closed on $\mathcal{D}$ with respect to
$T\pi_Q: TT^* Q \rightarrow TQ. $ Then one can construct
a magnetic symplectic form on $T^*Q$,
$\omega^B= \omega+ \pi_Q^*(\mathbf{d}\gamma), $ where $B=- \mathbf{d}\gamma,$
and a nonholonomic magnetic Hamiltonian system $(T^*Q, \omega^B, \mathcal{D},H )$
with an associated distributional magnetic Hamiltonian system
$(\mathcal{K},\omega^B_{\mathcal{K}},H_{\mathcal{K}})$.
Denote $X^\gamma = T\pi_{Q}\cdot X^B_H \cdot \gamma$,
where $X^B_H$ is the dynamical vector field
of the magnetic Hamiltonian system $(T^*Q,\omega^B,H)$.
Then $\gamma$ is a solution of the Type I of
Hamilton-Jacobi equation $T\gamma\cdot X^\gamma= X^B_{\mathcal{K}}\cdot \gamma ,$
for the distributional magnetic Hamiltonian system
$(\mathcal{K},\omega^B_{\mathcal {K}},H_{\mathcal{K}})$.
\end{theo}

Next, for any symplectic map $\varepsilon: T^* Q \rightarrow T^* Q $
with respect to the magnetic symplectic form $\omega^B$,
we can prove the following Type II of
Hamilton-Jacobi theorem for the distributional magnetic Hamiltonian system.
For convenience, the
maps involved in the following theorem and its proof are shown in
Diagram-4.

\begin{center}
\hskip 0cm \xymatrix{& \mathcal{M} \ar[d]_{X^B_{\mathcal{K}}}
\ar[r]^{i_{\mathcal{M}}} & T^* Q \ar[d]_{X^B_{H\cdot \varepsilon}}
\ar[dr]^{X^\varepsilon} \ar[r]^{\pi_Q}
& Q \ar[r]^{\gamma} & T^*Q \ar[d]^{X^B_H} \\
& \mathcal{K}  & T(T^*Q) \ar[l]_{\tau_{\mathcal{K}}} & TQ
\ar[l]_{T\gamma} & T(T^* Q)\ar[l]_{T\pi_Q}}
\end{center}
$$\mbox{Diagram-4}$$

\begin{theo} (Type II of Hamilton-Jacobi Theorem for a Distributional Magnetic Hamiltonian System)
For the nonholonomic magnetic Hamiltonian system $(T^*Q,\omega^B,\mathcal{D},H)$
with an associated distributional magnetic Hamiltonian system
$(\mathcal{K},\omega^B_{\mathcal {K}},H_{\mathcal{K}})$, assume that $\gamma: Q
\rightarrow T^*Q$ is an one-form on $Q$, and $\lambda=\gamma \cdot
\pi_{Q}: T^* Q \rightarrow T^* Q, $ and for any
symplectic map $\varepsilon: T^* Q \rightarrow T^* Q $ with respect to
the magnetic symplectic form $\omega^B$, denote by
$X^\varepsilon = T\pi_{Q}\cdot X^B_H \cdot \varepsilon$,
where $X^B_{H}$ is the dynamical
vector field of the magnetic Hamiltonian system $(T^*Q,\omega^B,H)$.
Moreover, assume that $\textmd{Im}(\gamma)\subset
\mathcal{M}=\mathcal{F}L(\mathcal{D}), $ and $\varepsilon(\mathcal{M})\subset \mathcal{M},$
and $\textmd{Im}(T\gamma)\subset \mathcal{K}. $
Then $\varepsilon$ is a solution of the equation
$\tau_{\mathcal{K}}\cdot T\varepsilon(X^B_{H\cdot\varepsilon})= T\lambda \cdot X^B_H\cdot\varepsilon,$
if and only if it is a solution of the equation
$T\gamma \cdot X^\varepsilon= X^B_{\mathcal{K}} \cdot \varepsilon $.
Here $ X^B_{H\cdot\varepsilon}$ is the magnetic Hamiltonian vector field of the function
$H \cdot \varepsilon: T^* Q\rightarrow \mathbb{R}, $ and $X^B_{\mathcal{K}}$
is the dynamical vector field of the distributional magnetic Hamiltonian system
$(\mathcal{K},\omega^B_{\mathcal {K}},H_{\mathcal{K}})$. The equation $T\gamma \cdot
X^\varepsilon= X^B_{\mathcal{K}} \cdot \varepsilon,$ is called the Type II of
Hamilton-Jacobi equation for the distributional magnetic Hamiltonian system
$(\mathcal{K},\omega^B_{\mathcal {K}},H_{\mathcal{K}})$.
\end{theo}

\noindent{\bf Proof: } In the same way, we note that
$\textmd{Im}(\gamma)\subset \mathcal{M}, $ and
$\textmd{Im}(T\gamma)\subset \mathcal{K}, $ in this case,
$\omega^B_{\mathcal{K}}\cdot
\tau_{\mathcal{K}}=\tau_{\mathcal{K}}\cdot \omega^B_{\mathcal{M}}=
\tau_{\mathcal{K}}\cdot i_{\mathcal{M}}^* \cdot\omega^B, $ along
$\textmd{Im}(T\gamma)$. Moreover, from the distributional magnetic Hamiltonian equation (2.2),
we have that $X^B_{\mathcal{K}}= \tau_{\mathcal{K}}\cdot X^B_H,$
and $\tau_{\mathcal{K}}\cdot X^B_{H}\cdot \varepsilon = X^B_{\mathcal{K}}\cdot \varepsilon $.
Note that $\varepsilon(\mathcal{M})\subset \mathcal{M},$ and
$T\pi_{Q}(X^B_H\cdot \varepsilon(q,p))\in
\mathcal{D}_{q}, \; \forall q \in Q, \; (q,p) \in \mathcal{M}(\subset T^* Q), $
and hence $X^B_H\cdot \varepsilon \in \mathcal{F}$ along $\varepsilon$.
Thus, using the non-degenerate distributional two-form
$\omega^B_{\mathcal{K}}$, from Lemma 3.4 and Lemma 4.3, if
we take that $v= X^B_H\cdot \varepsilon \in \mathcal{F}, $
and for any $w \in \mathcal{F}, \; T\lambda(w)\neq 0, $ and
$\tau_{\mathcal{K}}\cdot w \neq 0, $ then we have that
\begin{align*}
& \omega^B_{\mathcal{K}}(T\gamma \cdot X^\varepsilon, \;
\tau_{\mathcal{K}}\cdot w)=
\omega^B_{\mathcal{K}}(\tau_{\mathcal{K}}\cdot T\gamma \cdot
X^\varepsilon, \; \tau_{\mathcal{K}}\cdot w)\\ & =
\tau_{\mathcal{K}}\cdot i_{\mathcal{M}}^* \cdot\omega^B(T\gamma \cdot
T\pi_Q \cdot X^B_H \cdot \varepsilon, \; w ) = \tau_{\mathcal{K}}\cdot
i_{\mathcal{M}}^* \cdot\omega^B(T(\gamma \cdot \pi_Q)\cdot X^B_H \cdot \varepsilon, \; w)\\
& =\tau_{\mathcal{K}}\cdot i_{\mathcal{M}}^* \cdot(\omega^B(X^B_H \cdot
\varepsilon, \; w-T(\gamma \cdot \pi_Q)\cdot w)
-(\mathbf{d}\gamma+B)(T\pi_{Q}(X^B_H\cdot \varepsilon), \; T\pi_{Q}(w)))\\
& = \tau_{\mathcal{K}}\cdot i_{\mathcal{M}}^* \cdot\omega^B(X^B_H \cdot
\varepsilon, \; w) - \tau_{\mathcal{K}}\cdot i_{\mathcal{M}}^* \cdot
\omega^B(X^B_H \cdot \varepsilon, \; T\lambda \cdot w)\\
& \;\;\;\;\;\;
-\tau_{\mathcal{K}}\cdot i_{\mathcal{M}}^* \cdot (\mathbf{d}\gamma+B)
(T\pi_{Q}(X^B_H\cdot \varepsilon), \; T\pi_{Q}(w))\\
& = \omega^B_{\mathcal{K}}( \tau_{\mathcal{K}}\cdot X^B_H \cdot \varepsilon,
\; \tau_{\mathcal{K}}\cdot w) -
\omega^B_{\mathcal{K}}(\tau_{\mathcal{K}}\cdot X^B_H \cdot \varepsilon, \;
\tau_{\mathcal{K}}\cdot T\lambda \cdot w)\\
& \;\;\;\;\;\;
+\tau_{\mathcal{K}}\cdot i_{\mathcal{M}}^* \cdot \lambda^* \omega^B(X^B_H\cdot \varepsilon, \; w)\\
& = \omega^B_{\mathcal{K}}(X^B_{\mathcal{K}}\cdot \varepsilon, \;
\tau_{\mathcal{K}} \cdot w) -
\omega^B_{\mathcal{K}}(\tau_{\mathcal{K}}\cdot X^B_H \cdot \varepsilon,
\; T\lambda \cdot w)+ \omega^B_{\mathcal{K}}(T\lambda\cdot X^B_H\cdot \varepsilon,
\; T\lambda \cdot w),
\end{align*}
where we have used that $ \tau_{\mathcal{K}}\cdot T\gamma
= T\gamma, \; \tau_{\mathcal{K}}\cdot T\lambda= T\lambda, $ and
$\tau_{\mathcal{K}}\cdot X^B_H\cdot \varepsilon = X^B_{\mathcal{K}}\cdot
\varepsilon, $ since $\textmd{Im}(T\gamma)\subset \mathcal{K}. $ Note
that $\varepsilon: T^* Q \rightarrow T^* Q $ is symplectic
with respect to the magnetic symplectic form $\omega^B$, and $
X^B_H\cdot \varepsilon = T\varepsilon \cdot X^B_{H\cdot\varepsilon}, $ along
$\varepsilon$, and hence $\tau_{\mathcal{K}}\cdot X^B_H \cdot \varepsilon=
\tau_{\mathcal{K}}\cdot T\varepsilon \cdot X^B_{H \cdot \varepsilon}, $ along $\varepsilon$.
Then we have that
\begin{align*}
& \omega^B_{\mathcal{K}}(T\gamma \cdot X^\varepsilon, \;
\tau_{\mathcal{K}}\cdot w)-
\omega^B_{\mathcal{K}}(X^B_{\mathcal{K}}\cdot \varepsilon, \;
\tau_{\mathcal{K}} \cdot w) \nonumber \\
& = - \omega^B_{\mathcal{K}}(\tau_{\mathcal{K}}\cdot X^B_H \cdot \varepsilon, \;
 T\lambda \cdot w)+ \omega^B_{\mathcal{K}}(T\lambda\cdot X^B_H\cdot \varepsilon,
\; T\lambda \cdot w)\\
&= \omega^B_{\mathcal{K}}(T\lambda\cdot X^B_H\cdot \varepsilon
-\tau_{\mathcal{K}}\cdot T\varepsilon \cdot X^B_{H \cdot \varepsilon},
\; T\lambda \cdot w).
\end{align*}
Because the induced distributional two-form
$\omega^B_{\mathcal{K}}$ is non-degenerate, it follows that the equation
$T\gamma\cdot X^\varepsilon= X^B_{\mathcal{K}}\cdot
\varepsilon ,$ is equivalent to the equation
$\tau_{\mathcal{K}}\cdot T\varepsilon \cdot X^B_{H\cdot\varepsilon} = T\lambda\cdot X^B_H\cdot \varepsilon $.
Thus, $\varepsilon$ is a solution of the equation
$\tau_{\mathcal{K}}\cdot T\varepsilon\cdot X^B_{H\cdot\varepsilon}= T\lambda \cdot X^B_H \cdot\varepsilon,$
if and only if it is a solution of
the Type II of Hamilton-Jacobi equation $T\gamma\cdot X^\varepsilon= X^B_{\mathcal{K}}\cdot
\varepsilon .$
\hskip 0.3cm $\blacksquare$

\begin{rema}
It is worthy of noting that, the Type I of Hamilton-Jacobi equation
$T\gamma \cdot X^\gamma= X^B_{\mathcal{K}}\cdot \gamma ,$
is the equation of the differential one-form $\gamma$; and
the Type II of Hamilton-Jacobi equation $T\gamma \cdot X^\varepsilon
= X^B_{\mathcal{K}}\cdot \varepsilon ,$ is the equation of
the symplectic diffeomorphism map $\varepsilon$.
If the nonholonomic magnetic Hamiltonian system we considered has not any constrains, in this case,
the distributional magnetic Hamiltonian system is just the magnetic Hamiltonian system itself.
From the above Type I and Type II of Hamilton-Jacobi theorems, that is,
Theorem 4.4 and Theorem 4.6, we can get the Theorem 3.5 and Theorem 3.7.
It shows that Theorem 4.4 and Theorem 4.6 can be regarded as an extension of two types of
Hamilton-Jacobi theorem for the magnetic Hamiltonian system to the system with
nonholonomic context. On the other hand, when $B=0$,
in this case the magnetic symplectic form $\omega^B$
is just the canonical symplectic form $\omega$ on $T^*Q$, and
the distributional magnetic Hamiltonian system is just the distributional Hamiltonian system itself.
From the above Type I and Type II of Hamilton-Jacobi theorems, that is,
Theorem 4.4 and Theorem 4.6, we can get the Theorem 3.5 and Theorem 3.6
given by Le\'{o}n and Wang in \cite{lewa15}.
It shows that Theorem 4.4 and Theorem 4.6 can be regarded as an extension of two types of
Hamilton-Jacobi theorem for the distributional Hamiltonian system to that for the
distributional magnetic Hamiltonian system.
\end{rema}

\section{Nonholonomic Reduced Distributional Magnetic Hamiltonian System }

It is well-known that the reduction theory for the mechanical system with symmetry
is an important subject and it is widely studied in the theory of
mathematics and mechanics, as well as applications; see
Abraham and Marsden \cite{abma78}, Arnold \cite{ar89},
Libermann and Marle \cite{lima87}, Marsden \cite{ma92}, Marsden et al.
\cite{mamiorpera07, mamora90},
Marsden and Ratiu \cite{mara99}, Marsden and
Weinstein \cite{mawe74}, and Ortega and Ratiu \cite{orra04}
and so on, for more details and development.
In particular, the reduction of nonholonomically constrained mechanical systems
is also very important subject in geometric mechanics, and it is
regarded as a useful tool for simplifying and studying
concrete nonholonomic systems, see
Bates and $\acute{S}$niatycki \cite{basn93},  Cantrijn et al.
\cite{calemama99}, Cendra et al. \cite{cemara01},
Cushman et al. \cite{cudusn10} and \cite{cukesnba95}, Koiller \cite{ko92},
Le\'{o}n and Rodrigues \cite{lero89} and Le\'{o}n and
Wang \cite{lewa15} and so on.\\

In this section, we shall consider the nonholonomic reduction and
Hamilton-Jacobi theory of a nonholonomic magnetic Hamiltonian
system with symmetry. At first, we give the definition of
a nonholonomic magnetic Hamiltonian system with symmetry.
Then, by using the similar method in Le\'{o}n and Wang \cite{lewa15}
and Bates and $\acute{S}$niatycki \cite{basn93}.
and by analyzing carefully the dynamics and structure of
the nonholonomic magnetic Hamiltonian system with symmetry,
we give a geometric formulation of
the nonholonomic reduced distributional magnetic Hamiltonian system,
Moreover, we derive precisely the geometric constraint conditions of
the non-degenerate, and nonholonomic reduced distributional two-form
for the nonholonomic reducible dynamical vector field,
that is, the two types of Hamilton-Jacobi equation for the
nonholonomic reduced distributional magnetic Hamiltonian system,
which are an extension of the above two types of Hamilton-Jacobi equation
for the distributional magnetic Hamiltonian system
given in section 4 under nonholonomic reduction.\\

Assume that the Lie group $G$ acts smoothly on the manifold $Q$ by the left,
and we also consider the natural lifted actions on $TQ$ and $T^* Q$,
and assume that the cotangent lifted action on $T^\ast Q$ is free, proper and
symplectic with respect to the magnetic symplectic form
$\omega^B= \omega- \pi_Q^*B $ on $T^*Q$,
where $\omega$ is the canonical symplectic form on $T^* Q$
and $B$ is a closed two-form on $Q$.
The orbit space $T^* Q/ G$ is a smooth manifold and the
canonical projection $\pi_{/G}: T^* Q \rightarrow T^* Q /G $ is
a surjective submersion. Assume that $H: T^*Q \rightarrow \mathbb{R}$ is a
$G$-invariant Hamiltonian, and the $\mathcal{D}$-completely and
$\mathcal{D}$-regularly nonholonomic constraint $\mathcal{D}\subset
TQ$ is a $G$-invariant distribution, that is, the tangent of group action maps
$\mathcal{D}_q$ to $\mathcal{D}_{gq}$ for any
$q\in Q $. A nonholonomic magnetic Hamiltonian system with symmetry
is 5-tuple $(T^*Q,G,\omega^B,\mathcal{D},H)$, which is a
magnetic Hamiltonian system with symmetry and $G$-invariant
nonholonomic constraint $\mathcal{D}$.\\

In the following we first consider the nonholonomic reduction
of a nonholonomic magnetic Hamiltonian system with symmetry
$(T^*Q,G,\omega^B,\mathcal{D},H)$.
Note that the Legendre transformation $\mathcal{F}L: TQ
\rightarrow T^*Q$ is a fiber-preserving map,
and $\mathcal{D}\subset TQ$ is $G$-invariant
for the tangent lifted left action $\Phi^{T}: G\times TQ\rightarrow TQ, $
then the constraint submanifold
$\mathcal{M}=\mathcal{F}L(\mathcal{D})\subset T^*Q$ is
$G$-invariant for the cotangent lifted left action $\Phi^{T^\ast}:
G\times T^\ast Q\rightarrow T^\ast Q$.
For the nonholonomic magnetic Hamiltonian system with symmetry
$(T^*Q,G, \omega^B, \mathcal{D},H )$,
in the same way, we define the distribution $\mathcal{F}$, which is the pre-image of the
nonholonomic constraints $\mathcal{D}$ for the map $T\pi_Q: TT^* Q
\rightarrow TQ$, that is, $\mathcal{F}=(T\pi_Q)^{-1}(\mathcal{D})$,
and the distribution $\mathcal{K}=\mathcal{F} \cap T\mathcal{M}$.
Moreover, we can also define the distributional magnetic two-form $\omega^B_\mathcal{K}$,
which is induced from the magnetic symplectic form $\omega^B$ on $T^* Q$, that is,
$\omega^B_\mathcal{K}= \tau_{\mathcal{K}}\cdot \omega^B_{\mathcal{M}},$ and
$\omega^B_{\mathcal{M}}= i_{\mathcal{M}}^* \omega^B $.
If the admissibility condition $\mathrm{dim}\mathcal{M}=
\mathrm{rank}\mathcal{F}$ and the compatibility condition
$T\mathcal{M}\cap \mathcal{F}^\bot= \{0\}$ hold, then
$\omega^B_\mathcal{K}$ is non-degenerate as a
bilinear form on each fibre of $\mathcal{K}$, there exists a vector
field $X^B_\mathcal{K}$ on $\mathcal{M}$ which takes values in the
constraint distribution $\mathcal{K}$, such that for the function $H_\mathcal{K}$,
the following distributional magnetic Hamiltonian equation holds, that is,
\begin{align}
\mathbf{i}_{X^B_\mathcal{K}}\omega^B_\mathcal{K}
=\mathbf{d}H_\mathcal{K},
\label{5.1} \end{align}
where the function $H_{\mathcal{K}}$ satisfies
$\mathbf{d}H_{\mathcal{K}}= \tau_{\mathcal{K}}\cdot \mathbf{d}H_{\mathcal {M}}$,
and $H_\mathcal{M}= \tau_{\mathcal{M}}\cdot H$
is the restriction of $H$ to $\mathcal{M}$, and
from the equation (5.1), we have that
$X^B_{\mathcal{K}}=\tau_{\mathcal{K}}\cdot X^B_H $.\\

In the following we define that the quotient space
$\bar{\mathcal{M}}=\mathcal{M}/G$ of the $G$-orbit in $\mathcal{M}$
is a smooth manifold with projection $\pi_{/G}:
\mathcal{M}\rightarrow \bar{\mathcal{M}}( \subset T^* Q /G),$ which
is a surjective submersion. The reduced magnetic symplectic form
$\omega^B_{\bar{\mathcal{M}}}= \pi^*_{/G} \cdot \omega^B_{\mathcal{M}}$
on $\bar{\mathcal{M}}$ is induced from the magnetic symplectic form
$\omega^B_{\mathcal{M}}= i_{\mathcal{M}}^* \omega^B $ on $\mathcal{M}$.
Since $G$ is the symmetry group of the system
$(T^*Q,G,\omega^B,\mathcal{D},H)$, all intrinsically
defined vector fields and distributions are pushed down to
$\bar{\mathcal{M}}$. In particular, the vector field $X^B_\mathcal{M}$
on $\mathcal{M}$ is pushed down to a vector field
$X^B_{\bar{\mathcal{M}}}=T\pi_{/G}\cdot X^B_\mathcal{M}$, and the
distribution $\mathcal{K}$ is pushed down to a distribution
$T\pi_{/G}\cdot \mathcal{K}$ on $\bar{\mathcal{M}}$, and the
Hamiltonian $H$ is pushed down to $h_{\bar{\mathcal{M}}}$, such that
$h_{\bar{\mathcal{M}}}\cdot \pi_{/G}=
\tau_{\mathcal{M}}\cdot H$. However, $\omega^B_\mathcal{K}$ need not
to be pushed down to a distributional two-form defined on $T\pi_{/G}\cdot
\mathcal{K}$, despite of the fact that $\omega^B_\mathcal{K}$ is
$G$-invariant. This is because there may be infinitesimal symmetry
$\eta_{\mathcal{K}}$ that lies in $\mathcal{M}$, such that
$\mathbf{i}_{\eta_\mathcal{K}} \omega^B_\mathcal{K}\neq 0$. From Bates
and $\acute{S}$niatycki \cite{basn93}, we know that in order to eliminate
this difficulty, $\omega^B_\mathcal{K}$ is restricted to a
sub-distribution $\mathcal{U}$ of $\mathcal{K}$ defined by
$$\mathcal{U}=\{u\in\mathcal{K} \; | \; \omega^B_\mathcal{K}(u,v)
=0,\quad \forall \; v \in \mathcal{V}\cap \mathcal{K}\},$$ where
$\mathcal{V}$ is the distribution on $\mathcal{M}$ tangent to the
orbits of $G$ in $\mathcal{M}$ and it is spanned by the infinitesimal
symmetries. Clearly, $\mathcal{U}$ and $\mathcal{V}$ are both
$G$-invariant, project down to $\bar{\mathcal{M}}$ and
$T\pi_{/G}\cdot \mathcal{V}=0$, and define the distribution $\bar{\mathcal{K}}$ by
$\bar{\mathcal{K}}= T\pi_{/G}\cdot \mathcal{U}$. Moreover, we take
that $\omega^B_\mathcal{U}= \tau_{\mathcal{U}}\cdot
\omega^B_{\mathcal{M}}$ is the restriction of the induced magnetic symplectic form
$\omega^B_{\mathcal{M}}$ on $T^*\mathcal{M}$ fibrewise to the
distribution $\mathcal{U}$, where $\tau_{\mathcal{U}}$ is the
restriction map to distribution $\mathcal{U}$, and the
$\omega^B_{\mathcal{U}}$ is pushed down to a
distributional magnetic two-form $\omega^B_{\bar{\mathcal{K}}}$ on
$\bar{\mathcal{K}}$, such that $\pi_{/G}^*
\omega^B_{\bar{\mathcal{K}}}= \omega^B_{\mathcal{U}}$.
It is worthy of noting that the distributional magnetic two-form
$\omega^B_{\bar{\mathcal{K}}}$ is not a true two-form
on a manifold, so it does not make sense to speak about it being closed.
Thus, it is called the nonholonomic reduced
distributional magnetic two-form to avoid any confusion.\\

From the above construction we know that,
if the admissibility condition $\mathrm{dim}\bar{\mathcal{M}}=
\mathrm{rank}\bar{\mathcal{F}}$ and the compatibility condition
$T\bar{\mathcal{M}} \cap \bar{\mathcal{F}}^\bot= \{0\}$ hold, where
$\bar{\mathcal{F}}^\bot$ denotes the symplectic orthogonal of
$\bar{\mathcal{F}}$ with respect to the reduced magnetic symplectic form
$\omega^B_{\bar{\mathcal{M}}}$, then the nonholonomic reduced
distributional magnetic two-form
$\omega^B_{\bar{\mathcal{K}}}$ is non-degenerate as a bilinear form on
each fibre of $\bar{\mathcal{K}}$, and hence there exists a vector field
$X^B_{\bar{\mathcal{K}}}$ on $\bar{\mathcal{M}}$ which takes values in
the constraint distribution $\bar{\mathcal{K}}$, such that the nonholonomic
reduced distributional magnetic Hamiltonian equation holds, that is,
\begin{align}
\mathbf{i}_{X^B_{\bar{\mathcal{K}}}}\omega^B_{\bar{\mathcal{K}}}
=\mathbf{d}h_{\bar{\mathcal{K}}},
\label{5.2} \end{align}
where $\mathbf{d}h_{\bar{\mathcal{K}}}$ is the restriction of
$\mathbf{d}h_{\bar{\mathcal{M}}}$ to $\bar{\mathcal{K}}$ and
the function $h_{\bar{\mathcal{K}}}$ satisfies
$\mathbf{d}h_{\bar{\mathcal{K}}}
= \tau_{\bar{\mathcal{K}}}\cdot \mathbf{d}h_{\bar{\mathcal{M}}}$,
and $h_{\bar{\mathcal{M}}}\cdot \pi_{/G}= H_{\mathcal{M}}$ and
$H_{\mathcal{M}}$ is the restriction of the Hamiltonian function $H$
to $\mathcal{M}$. In addition, from the distributional magnetic Hamiltonian equation (2.2),
$\mathbf{i}_{X^B_\mathcal{K}}\omega^B_\mathcal{K}=\mathbf{d}H_\mathcal
{K},$ we have that $X^B_{\mathcal{K}}=\tau_{\mathcal{K}}\cdot X^B_H, $
and from the nonholonomic reduced distributional magnetic Hamiltonian equation (5.2),
$\mathbf{i}_{X^B_{\bar{\mathcal{K}}}}\omega^B_{\bar{\mathcal{K}}}
=\mathbf{d}h_{\bar{\mathcal{K}}}$, we have that
$X^B_{\bar{\mathcal{K}}}
=\tau_{\bar{\mathcal{K}}}\cdot X^B_{h_{\bar{\mathcal{K}}}},$
where $ X^B_{h_{\bar{\mathcal{K}}}}$ is the magnetic Hamiltonian vector field of
the function $h_{\bar{\mathcal{K}}}$ with respect to the reduced magnetic symplectic
form $\omega^B_{\bar{\mathcal{M}}}$,
and the vector fields $X^B_{\mathcal{K}}$
and $X^B_{\bar{\mathcal{K}}}$ are $\pi_{/G}$-related,
that is, $X^B_{\bar{\mathcal{K}}}\cdot \pi_{/G}=T\pi_{/G}\cdot X^B_{\mathcal{K}}.$
Thus, the geometrical formulation of a nonholonomic reduced distributional
magnetic Hamiltonian system may be summarized as follows.

\begin{defi} (Nonholonomic Reduced Distributional magnetic Hamiltonian System)
Assume that the 5-tuple $(T^*Q,G,\omega^B,\mathcal{D},H)$ is a nonholonomic
magnetic Hamiltonian system with symmetry, where $\omega^B$ is the magnetic
symplectic form on $T^* Q$, and $\mathcal{D}\subset TQ$ is a
$\mathcal{D}$-completely and $\mathcal{D}$-regularly nonholonomic
constraint of the system, and $\mathcal{D}$ and $H$ are both
$G$-invariant. If there exists a nonholonomic reduced distribution $\bar{\mathcal{K}}$,
an associated non-degenerate  and nonholonomic reduced
distributional two-form $\omega^B_{\bar{\mathcal{K}}}$
and a vector field $X^B_{\bar{\mathcal {K}}}$ on the reduced constraint
submanifold $\bar{\mathcal{M}}=\mathcal{M}/G, $ where
$\mathcal{M}=\mathcal{F}L(\mathcal{D})\subset T^*Q$, such that the
nonholonomic reduced distributional magnetic Hamiltonian equation
$\mathbf{i}_{X^B_{\bar{\mathcal{K}}}}\omega^B_{\bar{\mathcal{K}}} =
\mathbf{d}h_{\bar{\mathcal{K}}}$ holds,
where $\mathbf{d}h_{\bar{\mathcal{K}}}$ is the restriction of
$\mathbf{d}h_{\bar{\mathcal{M}}}$ to $\bar{\mathcal{K}}$ and
the function $h_{\bar{\mathcal{K}}}$ satisfies
$\mathbf{d}h_{\bar{\mathcal{K}}}= \tau_{\bar{\mathcal{K}}}\cdot \mathbf{d}h_{\bar{\mathcal{M}}}$
and $h_{\bar{\mathcal{M}}}\cdot \pi_{/G}= H_{\mathcal{M}}$ as defined above.
Then the triple $(\bar{\mathcal{K}},\omega^B_{\bar{\mathcal {K}}},h_{\bar{\mathcal{K}}})$
is called a nonholonomic reduced distributional magnetic Hamiltonian system
of the nonholonomic magnetic Hamiltonian system with symmetry
$(T^*Q,G,\omega^B,\mathcal{D},H)$, and $X^B_{\bar{\mathcal {K}}}$ is
called a nonholonomic reduced dynamical vector field.
of the system
$(\bar{\mathcal{K}},\omega^B_{\bar{\mathcal{K}}},h_{\bar{\mathcal{K}}})$. Under the above
circumstances, we refer to $(T^*Q,G,\omega^B,\mathcal{D},H)$ as a
nonholonomic reducible magnetic Hamiltonian system with the associated
distributional magnetic Hamiltonian system
$(\mathcal{K},\omega^B_{\mathcal {K}},H_{\mathcal{K}})$
and nonholonomic reduced distributional magnetic Hamiltonian system
$(\bar{\mathcal{K}},\omega^B_{\bar{\mathcal{K}}},h_{\bar{\mathcal{K}}})$.
\end{defi}

Since the non-degenerate and nonholonomic reduced distributional two-form
$\omega^B_{\bar{\mathcal{K}}}$ is not a "true two-form"
on a manifold, and it is not symplectic, and hence
the nonholonomic reduced distributional magnetic Hamiltonian system
$(\bar{\mathcal{K}},\omega^B_{\bar{\mathcal{K}}},h_{\bar{\mathcal{K}}})$
may not be yet a Hamiltonian system, and may have no generating function,
and hence we can not describe the Hamilton-Jacobi equation for the nonholonomic reduced
distributional magnetic Hamiltonian system just like as in Theorem 1.1.
But, since the nonholonomic reduced distributional magnetic Hamiltonian system is a
dynamical system closely related to a magnetic Hamiltonian system,
for a given nonholonomic reducible magnetic Hamiltonian system
$(T^*Q,G,\omega^B,\mathcal{D},H)$ with the associated
distributional magnetic Hamiltonian system
$(\mathcal{K},\omega^B_{\mathcal {K}},H_{\mathcal{K}})$
and the nonholonomic reduced distributional magnetic Hamiltonian system
$(\bar{\mathcal{K}},\omega^B_{\bar{\mathcal {K}}},h_{\bar{\mathcal{K}}})$,
by using Lemma 3.4 and Lemma 4.3,
we can also derive precisely
the geometric constraint conditions of the nonholonomic reduced distributional two-form
$\omega^B_{\bar{\mathcal{K}}}$ for the dynamical vector field $X^B_{\bar{\mathcal {K}}}$,
that is, the two types of Hamilton-Jacobi equation for the
nonholonomic reduced distributional magnetic Hamiltonian system
$(\bar{\mathcal{K}},\omega^B_{\bar{\mathcal {K}}},h_{\bar{\mathcal{K}}})$.
At first, using the fact that the one-form $\gamma: Q
\rightarrow T^*Q $ satisfies the condition that
$\mathbf{d}\gamma=-B$ on $\mathcal{D}$ with respect to
$T\pi_Q: TT^* Q \rightarrow TQ, $
$\textmd{Im}(\gamma)\subset \mathcal{M}, $ and it is $G$-invariant,
as well as $ \textmd{Im}(T\gamma)\subset \mathcal{K}, $
we can prove the Type I of
Hamilton-Jacobi theorem for the nonholonomic reduced distributional
magnetic Hamiltonian system. For convenience, the maps involved in the
following theorem and its proof are shown in Diagram-5.
\begin{center}
\hskip 0cm \xymatrix{ & \mathcal{M} \ar[d]_{X^B_{\mathcal{K}}}
\ar[r]^{i_{\mathcal{M}}} & T^* Q \ar[d]_{X^B_{H}}
 \ar[r]^{\pi_Q}
  & Q \ar[d]_{X^\gamma} \ar[r]^{\gamma}
  & T^* Q \ar[d]_{X^B_H} \ar[r]^{\pi_{/G}} & T^* Q/G \ar[d]_{X^B_{h_{\bar{\mathcal{M}}}}}
  & \mathcal{\bar{M}} \ar[l]_{i_{\mathcal{\bar{M}}}} \ar[d]_{X^B_{\mathcal{\bar{K}}}}\\
  & \mathcal{K}
  & T(T^*Q) \ar[l]_{\tau_{\mathcal{K}}}
  & TQ \ar[l]_{T\gamma}
  & T(T^* Q) \ar[l]_{T\pi_Q} \ar[r]^{T\pi_{/G}}
  & T(T^* Q/G) \ar[r]^{\tau_{\mathcal{\bar{K}}}} & \mathcal{\bar{K}} }
\end{center}
$$\mbox{Diagram-5}$$

\begin{theo} (Type I of Hamilton-Jacobi Theorem for a Nonholonomic
Reduced Distributional magnetic Hamiltonian System)
For a given nonholonomic reducible magnetic Hamiltonian system
$(T^*Q,G,\omega^B,\mathcal{D},H)$ with the associated
distributional magnetic Hamiltonian system
$(\mathcal{K},\omega^B_{\mathcal {K}},H_{\mathcal{K}})$
and the nonholonomic reduced distributional magnetic Hamiltonian system
$(\bar{\mathcal{K}},\omega^B_{\bar{\mathcal{K}}},h_{\bar{\mathcal{K}}})$, assume that
$\gamma: Q \rightarrow T^*Q$ is an one-form on $Q$, and
$X^\gamma = T\pi_{Q}\cdot X^B_H \cdot \gamma$, where $X^B_{H}$ is
the magnetic Hamiltonian vector field of the corresponding unconstrained
magnetic Hamiltonian system with symmetry $(T^*Q,G,\omega^B,H)$. Moreover,
assume that $\textmd{Im}(\gamma)\subset \mathcal{M}, $ and it is
$G$-invariant, $ \textmd{Im}(T\gamma)\subset \mathcal{K}, $ and
$\bar{\gamma}=\pi_{/G}(\gamma): Q \rightarrow T^* Q/G .$ If the
one-form $\gamma: Q \rightarrow T^*Q $ satisfies the condition that
$\mathbf{d}\gamma=-B$ on $\mathcal{D}$ with respect to
$T\pi_Q: TT^* Q \rightarrow TQ, $ then $\bar{\gamma}$ is a solution
of the equation $T\bar{\gamma}\cdot X^ \gamma =
X^B_{\bar{\mathcal{K}}}\cdot \bar{\gamma}. $
Here $X^B_{\bar{\mathcal{K}}}$ is the dynamical vector
field of the nonholonomic reduced distributional magnetic Hamiltonian system
$(\bar{\mathcal{K}},\omega^B_{\bar{\mathcal{K}}},h_{\bar{\mathcal{K}}})$.
The equation $ T\bar{\gamma}\cdot X^ \gamma = X^B_{\bar{\mathcal{K}}}\cdot
\bar{\gamma},$ is called the Type I of Hamilton-Jacobi equation for
the nonholonomic reduced distributional magnetic Hamiltonian system
$(\bar{\mathcal{K}},\omega^B_{\bar{\mathcal{K}}},h_{\bar{\mathcal{K}}})$.
\end{theo}

\noindent{\bf Proof: } At first, from Theorem 4.4, we know that
$\gamma$ is a solution of the Hamilton-Jacobi equation
$T\gamma\cdot X^\gamma= X^B_{\mathcal{K}}\cdot \gamma .$ Next, we note that
$\textmd{Im}(\gamma)\subset \mathcal{M}, $ and it is $G$-invariant,
$ \textmd{Im}(T\gamma)\subset \mathcal{K}, $ and hence
$\textmd{Im}(T\bar{\gamma})\subset \bar{\mathcal{K}}, $ in this case,
$\pi^*_{/G}\cdot\omega^B_{\bar{\mathcal{K}}}\cdot\tau_{\bar{\mathcal{K}}}= \tau_{\mathcal{U}}\cdot
\omega^B_{\mathcal{M}}= \tau_{\mathcal{U}}\cdot i_{\mathcal{M}}^*\cdot
\omega^B, $ along $\textmd{Im}(T\bar{\gamma})$.
From the distributional magnetic Hamiltonian equation (2.2),
we have that $X^B_{\mathcal{K}}= \tau_{\mathcal{K}}\cdot X^B_H,$
and $\tau_{\mathcal{K}}\cdot X^B_{H}\cdot \gamma = X^B_{\mathcal{K}}\cdot \gamma $.
Because the vector fields $X^B_{\mathcal{K}}$
and $X^B_{\bar{\mathcal{K}}}$ are $\pi_{/G}$-related,
$T\pi_{/G}(X^B_{\mathcal{K}})=X^B_{\bar{\mathcal{K}}}\cdot \pi_{/G}$,
and hence $\tau_{\bar{\mathcal{K}}}\cdot T\pi_{/G}(X^B_{\mathcal{K}}\cdot \gamma)
=\tau_{\bar{\mathcal{K}}}\cdot (T\pi_{/G}(X^B_{\mathcal{K}}))\cdot (\gamma)
= \tau_{\bar{\mathcal{K}}}\cdot (X^B_{\bar{\mathcal{K}}}\cdot \pi_{/G})\cdot (\gamma)
= \tau_{\bar{\mathcal{K}}}\cdot X^B_{\bar{\mathcal{K}}}\cdot \pi_{/G}(\gamma)
= X^B_{\bar{\mathcal{K}}}\cdot \bar{\gamma}.$
Thus, using the non-degenerate, nonholonomic reduced
distributional two-form $\omega^B_{\bar{\mathcal{K}}}$,
from Lemma 3.4(ii) and Lemma 4.3, if we take that
$v= X^B_{H}\cdot \gamma \in \mathcal{F},$
and for any $w \in \mathcal{F}, \; T\lambda(w)\neq 0, $ and
$\tau_{\bar{\mathcal{K}}}\cdot T\pi_{/G}\cdot w \neq 0, $ then we have that
\begin{align*}
& \omega^B_{\bar{\mathcal{K}}}(T\bar{\gamma} \cdot X^\gamma, \;
\tau_{\bar{\mathcal{K}}}\cdot T\pi_{/G} \cdot w)
= \omega^B_{\bar{\mathcal{K}}}(\tau_{\bar{\mathcal{K}}}\cdot T(\pi_{/G} \cdot
\gamma) \cdot X^\gamma, \; \tau_{\bar{\mathcal{K}}}\cdot T\pi_{/G} \cdot w )\\
& = \pi^*_{/G}\cdot \omega^B_{\bar{\mathcal{K}}}\cdot\tau_{\bar{\mathcal{K}}}
(T\gamma \cdot X^\gamma, \; w)
= \tau_{\mathcal{U}}\cdot i^*_{\mathcal{M}} \cdot \omega^B (T\gamma \cdot
T\pi_Q \cdot X^B_H \cdot \gamma, \; w)\\
& = \tau_{\mathcal{U}}\cdot i^*_{\mathcal{M}} \cdot
\omega^B (T(\gamma \cdot \pi_Q)\cdot X^B_H \cdot \gamma, \; w) \\
& = \tau_{\mathcal{U}}\cdot i^*_{\mathcal{M}} \cdot
(\omega^B (X^B_H \cdot \gamma, \; w-T(\gamma \cdot \pi_Q)\cdot w)
- (\mathbf{d}\gamma+B)(T\pi_{Q}(X^B_H\cdot \gamma), \; T\pi_{Q}(w)))\\
& = \tau_{\mathcal{U}}\cdot i^*_{\mathcal{M}} \cdot \omega^B (X^B_H \cdot
\gamma, \; w) - \tau_{\mathcal{U}}\cdot i^*_{\mathcal{M}} \cdot
\omega^B (X^B_H \cdot \gamma, \; T(\gamma \cdot \pi_Q) \cdot w)\\
& \;\;\;\;\;\;
- \tau_{\mathcal{U}}\cdot i^*_{\mathcal{M}} \cdot(\mathbf{d}\gamma+B)
(T\pi_{Q}(X^B_H\cdot \gamma), \; T\pi_{Q}(w))\\
& =\pi^*_{/G}\cdot \omega^B_{\bar{\mathcal{K}}}\cdot\tau_{\bar{\mathcal{K}}}(X^B_H \cdot \gamma, \;
w) - \pi^*_{/G}\cdot \omega^B_{\bar{\mathcal{K}}}\cdot\tau_{\bar{\mathcal{K}}}(X^B_H \cdot \gamma,
\; T(\gamma \cdot \pi_Q) \cdot w)\\
& \;\;\;\;\;\; - \tau_{\mathcal{U}}\cdot i^*_{\mathcal{M}}
\cdot(\mathbf{d}\gamma+B)(T\pi_{Q}(X^B_H\cdot \gamma), \; T\pi_{Q}(w))\\
& = \omega^B_{\bar{\mathcal{K}}}(\tau_{\bar{\mathcal{K}}}\cdot T\pi_{/G}(X^B_H \cdot \gamma), \;
\tau_{\bar{\mathcal{K}}}\cdot T\pi_{/G} \cdot w)
- \omega^B_{\bar{\mathcal{K}}}(\tau_{\bar{\mathcal{K}}}\cdot T\pi_{/G}(X^B_H \cdot
\gamma), \; \tau_{\bar{\mathcal{K}}}\cdot T(\pi_{/G} \cdot\gamma) \cdot T\pi_{Q}(w))\\
& \;\;\;\;\;\; - \tau_{\mathcal{U}}\cdot i^*_{\mathcal{M}}
\cdot(\mathbf{d}\gamma+B)(T\pi_{Q}(X^B_H\cdot \gamma), \; T\pi_{Q}(w))\\
& = \omega^B_{\bar{\mathcal{K}}}(\tau_{\bar{\mathcal{K}}}\cdot T\pi_{/G}(X^B_H)\cdot
\pi_{/G}(\gamma), \; \tau_{\bar{\mathcal{K}}}\cdot T\pi_{/G} \cdot w) -
\omega^B_{\bar{\mathcal{K}}}(\tau_{\bar{\mathcal{K}}}\cdot T\pi_{/G}(X^B_H)\cdot \pi_{/G}(\gamma), \;
\tau_{\bar{\mathcal{K}}}\cdot T\bar{\gamma} \cdot T\pi_{Q}(w))\\
& \;\;\;\;\;\; - \tau_{\mathcal{U}}\cdot i^*_{\mathcal{M}}
\cdot(\mathbf{d}\gamma+B)(T\pi_{Q}(X^B_H\cdot \gamma), \; T\pi_{Q}(w))\\
& = \omega^B_{\bar{\mathcal{K}}}(X^B_{\bar{\mathcal{K}}} \cdot
\bar{\gamma}, \; \tau_{\bar{\mathcal{K}}}\cdot T\pi_{/G} \cdot w)-
\omega^B_{\bar{\mathcal{K}}}(X^B_{\bar{\mathcal{K}}} \cdot
\bar{\gamma}, \; T\bar{\gamma} \cdot T\pi_{Q}(w)) \\
& \;\;\;\;\;\; - \tau_{\mathcal{U}}\cdot
i^*_{\mathcal{M}} \cdot(\mathbf{d}\gamma+B)(T\pi_{Q}(X^B_H\cdot \gamma),
\; T\pi_{Q}(w)),
\end{align*}
where we have used that $\tau_{\bar{\mathcal{K}}}\cdot T\pi_{/G}(X^B_H\cdot \gamma)
=\tau_{\bar{\mathcal{K}}}\cdot X^B_{\bar{\mathcal{K}}}\cdot \bar{\gamma}=
X^B_{\bar{\mathcal{K}}}\cdot \bar{\gamma}, $ and
$\tau_{\bar{\mathcal{K}}}\cdot T\bar{\gamma}=T\bar{\gamma}, $ since
$\textmd{Im}(T\bar{\gamma})\subset \bar{\mathcal{K}}. $
Note that $X^B_{H}\cdot \gamma, \; w \in \mathcal{F},$ and
$T\pi_{Q}(X^B_H\cdot \gamma), \; T\pi_{Q}(w) \in \mathcal{D}.$
If the one-form $\gamma: Q \rightarrow T^*Q $ satisfies the condition that
$\mathbf{d}\gamma=-B$ on $\mathcal{D}$ with respect to
$T\pi_Q: TT^* Q \rightarrow TQ, $ then we have that
$(\mathbf{d}\gamma+B)(T\pi_{Q}(X^B_H\cdot \gamma), \; T\pi_{Q}(w))=0, $
and hence
$$
\tau_{\mathcal{U}}\cdot i_{\mathcal{M}}^* \cdot(\mathbf{d}\gamma+B)
(T\pi_{Q}(X^B_H\cdot \gamma), \; T\pi_{Q}(w))=0,
$$
Thus, we have that
\begin{equation}
\omega^B_{\bar{\mathcal{K}}}(T\bar{\gamma} \cdot X^\gamma, \;
\tau_{\bar{\mathcal{K}}}\cdot T\pi_{/G} \cdot w)
- \omega^B_{\bar{\mathcal{K}}}(X^B_{\bar{\mathcal{K}}} \cdot
\bar{\gamma}, \; \tau_{\bar{\mathcal{K}}}\cdot T\pi_{/G} \cdot w)
= -\omega^B_{\bar{\mathcal{K}}}(X^B_{\bar{\mathcal{K}}} \cdot
\bar{\gamma}, \; T\bar{\gamma} \cdot T\pi_{Q}(w)).
\label{5.3} \end{equation}
If $\bar{\gamma}$ satisfies the equation $
T\bar{\gamma}\cdot X^ \gamma = X^B_{\bar{\mathcal{K}}}\cdot
\bar{\gamma} ,$
from Lemma 3.4(i) we know that the right side of (5.3) becomes that
\begin{align*}
 -\omega^B_{\bar{\mathcal{K}}}(X^B_{\bar{\mathcal{K}}} \cdot
\bar{\gamma}, \; T\bar{\gamma} \cdot T\pi_{Q}(w))
& = -\omega^B_{\bar{\mathcal{K}}}\cdot\tau_{\bar{\mathcal{K}}}
(T\bar{\gamma}\cdot X^\gamma, \; T\bar{\gamma} \cdot T\pi_{Q}(w))\\
& = -\bar{\gamma}^*\omega^B_{\bar{\mathcal{K}}}\cdot\tau_{\bar{\mathcal{K}}}
(T\pi_{Q} \cdot X^B_{H} \cdot \gamma, \; T\pi_{Q}(w))\\
& = - \gamma^* \cdot \pi^*_{/G}\cdot \omega^B_{\bar{\mathcal{K}}}\cdot
\tau_{\bar{\mathcal{K}}}(T\pi_{Q} \cdot X^B_{H} \cdot \gamma, \; T\pi_{Q}(w))\\
& = - \gamma^* \cdot \tau_{\mathcal{U}}\cdot
i_{\mathcal{M}}^* \cdot \omega^B (T\pi_{Q}(X^B_{H}\cdot\gamma), \; T\pi_{Q}(w))\\
& = -\tau_{\mathcal{U}}\cdot i_{\mathcal{M}}^* \cdot\gamma^*
\omega^B ( T\pi_{Q}(X^B_{H}\cdot\gamma), \; T\pi_{Q}(w))\\
& = \tau_{\mathcal{U}}\cdot i_{\mathcal{M}}^* \cdot
(\mathbf{d}\gamma+B)(T\pi_{Q}( X^B_{H}\cdot\gamma ), \; T\pi_{Q}(w))=0,
\end{align*}
where $\gamma^*\cdot \tau_{\mathcal{U}}\cdot i^*_{\mathcal{M}}
\cdot \omega^B= \tau_{\mathcal{U}}\cdot i^*_{\mathcal{M}}
\cdot\gamma^*\cdot \omega^B, $
because $\textmd{Im}(\gamma)\subset \mathcal{M}. $
But, since the nonholonomic reduced distributional two-form
$\omega^B_{\bar{\mathcal{K}}}$ is non-degenerate,
the left side of (5.3) equals zero, only when
$\bar{\gamma}$ satisfies the equation $
T\bar{\gamma}\cdot X^ \gamma = X^B_{\bar{\mathcal{K}}}\cdot
\bar{\gamma} .$ Thus,
if the one-form $\gamma: Q \rightarrow T^*Q $ satisfies the condition that
$\mathbf{d}\gamma=-B$ on $\mathcal{D}$ with respect to
$T\pi_Q: TT^* Q \rightarrow TQ, $ then $\bar{\gamma}$
must be a solution of the Type I of Hamilton-Jacobi equation
$T\bar{\gamma}\cdot X^ \gamma = X^B_{\bar{\mathcal{K}}}\cdot
\bar{\gamma}. $
\hskip 0.3cm $\blacksquare$\\

Next, for any $G$-invariant symplectic map $\varepsilon: T^* Q \rightarrow T^* Q $
with respect to $\omega^B$, we can prove
the following Type II of Hamilton-Jacobi theorem
for the nonholonomic reduced distributional magnetic Hamiltonian system.
For convenience, the maps involved in the following theorem and its
proof are shown in Diagram-6.
\begin{center}
\hskip 0cm \xymatrix{ & \mathcal{M} \ar[d]_{X^B_{\mathcal{K}}}
\ar[r]^{i_{\mathcal{M}}} & T^* Q \ar[d]_{X^B_{H\cdot \varepsilon}}
\ar[dr]^{X^\varepsilon} \ar[r]^{\pi_Q}
  & Q \ar[r]^{\gamma} & T^* Q \ar[d]_{X^B_H} \ar[r]^{\pi_{/G}} & T^* Q/G \ar[d]_{X^B_{h_{\bar{\mathcal{M}}}}}
  & \mathcal{\bar{M}} \ar[l]_{i_{\mathcal{\bar{M}}}} \ar[d]_{X^B_{\mathcal{\bar{K}}}}\\
  & \mathcal{K}
  & T(T^*Q) \ar[l]_{\tau_{\mathcal{K}}}
  & TQ \ar[l]_{T\gamma}
  & T(T^* Q) \ar[l]_{T\pi_Q} \ar[r]^{T\pi_{/G}}
  & T(T^* Q/G) \ar[r]^{\tau_{\mathcal{\bar{K}}}} & \mathcal{\bar{K}} }
\end{center}
$$\mbox{Diagram-6}$$

\begin{theo} (Type II of Hamilton-Jacobi Theorem for a Nonholonomic
Reduced Distributional magnetic Hamiltonian System)
For a given nonholonomic reducible magnetic Hamiltonian system
$(T^*Q,G,\omega^B,\mathcal{D},H)$ with the associated
distributional magnetic Hamiltonian system
$(\mathcal{K},\omega^B_{\mathcal {K}},H_{\mathcal{K}})$
and the nonholonomic reduced distributional magnetic Hamiltonian system
$(\bar{\mathcal{K}},\omega^B_{\bar{\mathcal{K}}},h_{\bar{\mathcal{K}}})$, assume that
$\gamma: Q \rightarrow T^*Q$ is an one-form on $Q$, and $\lambda=
\gamma \cdot \pi_{Q}: T^* Q \rightarrow T^* Q, $ and for any $G$-invariant
symplectic map $\varepsilon: T^* Q \rightarrow T^* Q $ with respect to $\omega^B$, denote by
$X^\varepsilon = T\pi_{Q}\cdot X^B_H \cdot \varepsilon$, where $X^B_{H}$ is
the magnetic Hamiltonian vector field of the corresponding unconstrained
magnetic Hamiltonian system with symmetry $(T^*Q,G,\omega^B,H)$. Moreover,
assume that $\textmd{Im}(\gamma)\subset \mathcal{M}, $ and it is
$G$-invariant, $\varepsilon(\mathcal{M})\subset \mathcal{M}$,
$ \textmd{Im}(T\gamma)\subset \mathcal{K}, $ and
$\bar{\gamma}=\pi_{/G}(\gamma): Q \rightarrow T^* Q/G $, and
$\bar{\lambda}=\pi_{/G}(\lambda): T^* Q \rightarrow T^* Q/G, $ and
$\bar{\varepsilon}=\pi_{/G}(\varepsilon): T^* Q \rightarrow T^* Q/G. $ Then
$\varepsilon$ and $\bar{\varepsilon}$ satisfy the equation
$\tau_{\bar{\mathcal{K}}}\cdot T\bar{\varepsilon}\cdot X^B_{h_{\bar{\mathcal{K}}}\cdot
\bar{\varepsilon}}= T\bar{\lambda} \cdot X^B_H\cdot \varepsilon, $ if and only if they satisfy the
equation $T\bar{\gamma}\cdot X^ \varepsilon =
X^B_{\bar{\mathcal{K}}}\cdot \bar{\varepsilon}. $ Here
$ X^B_{h_{\bar{\mathcal{K}}} \cdot\bar{\varepsilon}}$ is the magnetic Hamiltonian
vector field of the function $h_{\bar{\mathcal{K}}}\cdot \bar{\varepsilon}: T^* Q\rightarrow
\mathbb{R}, $ and $X^B_{\bar{\mathcal{K}}}$ is the dynamical vector
field of the nonholonomic reduced distributional magnetic Hamiltonian system
$(\bar{\mathcal{K}},\omega^B_{\bar{\mathcal{K}}},h_{\bar{\mathcal{K}}})$.
The equation $ T\bar{\gamma}\cdot X^\varepsilon
= X^B_{\bar{\mathcal{K}}}\cdot \bar{\varepsilon},$ is called
the Type II of Hamilton-Jacobi equation for the
nonholonomic reduced distributional magnetic Hamiltonian system
$(\bar{\mathcal{K}},\omega^B_{\bar{\mathcal{K}}},h_{\bar{\mathcal{K}}})$.
\end{theo}

\noindent{\bf Proof: } In the same way, we note that
$\textmd{Im}(\gamma)\subset \mathcal{M}, $ and it is $G$-invariant,
$ \textmd{Im}(T\gamma)\subset \mathcal{K}, $ and hence
$\textmd{Im}(T\bar{\gamma})\subset \bar{\mathcal{K}}, $ in this case,
$\pi^*_{/G}\cdot\omega^B_{\bar{\mathcal{K}}}\cdot \tau_{\bar{\mathcal{K}}}
= \tau_{\mathcal{U}}\cdot \omega^B_{\mathcal{M}}
= \tau_{\mathcal{U}}\cdot i_{\mathcal{M}}^*\cdot \omega^B, $
along $\textmd{Im}(T\bar{\gamma})$.
Moreover, from the distributional magnetic Hamiltonian equation (2.2),
we have that $X^B_{\mathcal{K}}= \tau_{\mathcal{K}}\cdot X^B_H.$
Note that $\varepsilon(\mathcal{M})\subset \mathcal{M},$ and
$T\pi_{Q}(X^B_H\cdot \varepsilon(q,p))\in
\mathcal{D}_{q}, \; \forall q \in Q, \; (q,p) \in \mathcal{M}(\subset T^* Q), $
and hence $X^B_H\cdot \varepsilon \in \mathcal{F}$ along $\varepsilon$.
Because the vector fields $X^B_{\mathcal{K}}$
and $X^B_{\bar{\mathcal{K}}}$ are $\pi_{/G}$-related,
$T\pi_{/G}(X^B_{\mathcal{K}})=X^B_{\bar{\mathcal{K}}}\cdot \pi_{/G}$,
and hence $\tau_{\bar{\mathcal{K}}}\cdot T\pi_{/G}(X^B_{\mathcal{K}}\cdot \varepsilon)
=\tau_{\bar{\mathcal{K}}}\cdot (T\pi_{/G}(X^B_{\mathcal{K}}))\cdot (\varepsilon)
= \tau_{\bar{\mathcal{K}}}\cdot (X^B_{\bar{\mathcal{K}}}\cdot \pi_{/G})\cdot (\varepsilon)
= \tau_{\bar{\mathcal{K}}}\cdot X^B_{\bar{\mathcal{K}}}\cdot \pi_{/G}(\varepsilon)
= X^B_{\bar{\mathcal{K}}}\cdot \bar{\varepsilon}.$ Thus, using the
non-degenerate and nonholonomic reduced distributional two-form $\omega^B_{\bar{\mathcal{K}}}$,
from Lemma 3.4 and Lemma 4.3, if we take that
$v=X^B_H\cdot \varepsilon \in \mathcal{F},$
and for any $w \in \mathcal{F}, \; T\lambda(w)\neq 0, $ and
$\tau_{\bar{\mathcal{K}}}\cdot T\pi_{/G}\cdot w \neq 0, $ then we have that
\begin{align*}
& \omega^B_{\bar{\mathcal{K}}}(T\bar{\gamma} \cdot X^\varepsilon, \;
\tau_{\bar{\mathcal{K}}}\cdot T\pi_{/G} \cdot w)
= \omega^B_{\bar{\mathcal{K}}}(\tau_{\bar{\mathcal{K}}}\cdot T(\pi_{/G} \cdot
\gamma) \cdot X^\varepsilon, \; \tau_{\bar{\mathcal{K}}}\cdot T\pi_{/G} \cdot w )\\
& = \pi^*_{/G}\cdot \omega^B_{\bar{\mathcal{K}}}
\cdot\tau_{\bar{\mathcal{K}}}(T\gamma \cdot X^\varepsilon, \; w) =
\tau_{\mathcal{U}}\cdot i^*_{\mathcal{M}} \cdot\omega^B (T\gamma \cdot
X^\varepsilon, \; w)\\
& = \tau_{\mathcal{U}}\cdot i^*_{\mathcal{M}} \cdot
\omega^B (T(\gamma \cdot \pi_Q)\cdot X^B_H \cdot \varepsilon, \; w) \\
& = \tau_{\mathcal{U}}\cdot i^*_{\mathcal{M}} \cdot
(\omega^B (X^B_H \cdot \varepsilon, \; w-T(\gamma \cdot \pi_Q)\cdot w)
- (\mathbf{d}\gamma+B)(T\pi_{Q}(X^B_H\cdot \varepsilon), \; T\pi_{Q}(w)))\\
& = \tau_{\mathcal{U}}\cdot i^*_{\mathcal{M}} \cdot \omega^B(X^B_H \cdot
\varepsilon, \; w) - \tau_{\mathcal{U}}\cdot i^*_{\mathcal{M}} \cdot
\omega^B(X^B_H \cdot \varepsilon, \; T\lambda \cdot w)\\
& \;\;\;\;\;\; - \tau_{\mathcal{U}}\cdot i^*_{\mathcal{M}} \cdot
(\mathbf{d}\gamma+B)(T\pi_{Q}(X^B_H\cdot \varepsilon), \; T\pi_{Q}(w))\\
& =\pi^*_{/G}\cdot \omega^B_{\bar{\mathcal{K}}}\cdot \tau_{\bar{\mathcal{K}}}(X^B_H \cdot \varepsilon, \;
w) - \pi^*_{/G}\cdot \omega^B_{\bar{\mathcal{K}}}\cdot \tau_{\bar{\mathcal{K}}}(X^B_H \cdot \varepsilon,
\; T\lambda \cdot w)+ \tau_{\mathcal{U}}\cdot i^*_{\mathcal{M}}
\cdot \lambda^* \omega^B(X^B_H\cdot \varepsilon, \; w)\\
& = \omega^B_{\bar{\mathcal{K}}}(\tau_{\bar{\mathcal{K}}}\cdot T\pi_{/G}(X^B_H \cdot \varepsilon), \;
\tau_{\bar{\mathcal{K}}}\cdot T\pi_{/G} \cdot w)
- \omega^B_{\bar{\mathcal{K}}}(\tau_{\bar{\mathcal{K}}}\cdot T\pi_{/G}(X^B_H \cdot
\varepsilon), \; \tau_{\bar{\mathcal{K}}}\cdot T(\pi_{/G} \cdot\lambda) \cdot w)\\
& \;\;\;\;\;\; +\pi^*_{/G}\cdot \omega^B_{\bar{\mathcal{K}}}\cdot \tau_{\bar{\mathcal{K}}}
(T\lambda\cdot X^B_H \cdot \varepsilon, \; T\lambda \cdot w)\\
& = \omega^B_{\bar{\mathcal{K}}}(\tau_{\bar{\mathcal{K}}}\cdot T\pi_{/G}(X^B_H)\cdot
\pi_{/G}(\varepsilon), \; \tau_{\bar{\mathcal{K}}}\cdot T\pi_{/G} \cdot w) -
\omega^B_{\bar{\mathcal{K}}}(\tau_{\bar{\mathcal{K}}}\cdot T\pi_{/G}(X^B_H)\cdot \pi_{/G}(\varepsilon), \;
\tau_{\bar{\mathcal{K}}}\cdot T\bar{\lambda} \cdot w)\\
& \;\;\;\;\;\; +\omega^B_{\bar{\mathcal{K}}}(\tau_{\bar{\mathcal{K}}}
\cdot T\pi_{/G}\cdot T\lambda\cdot X^B_H \cdot \varepsilon,
\; \tau_{\bar{\mathcal{K}}}\cdot T\pi_{/G}\cdot T\lambda \cdot w)\\
& = \omega^B_{\bar{\mathcal{K}}}(X^B_{\bar{\mathcal{K}}} \cdot
\bar{\varepsilon}, \; \tau_{\bar{\mathcal{K}}}\cdot T\pi_{/G} \cdot w)-
\omega^B_{\bar{\mathcal{K}}}(X^B_{\bar{\mathcal{K}}}\cdot
\bar{\varepsilon}, \; T\bar{\lambda} \cdot w)
+ \omega^B_{\bar{\mathcal{K}}}(T\bar{\lambda}\cdot X^B_H \cdot \varepsilon,
\; T\bar{\lambda} \cdot w),
\end{align*}
where we have used that $\tau_{\bar{\mathcal{K}}}\cdot T\pi_{/G}(X^B_H \cdot \varepsilon)
=\tau_{\bar{\mathcal{K}}}(X^B_{\bar{\mathcal{K}}})\cdot \bar{\varepsilon}=
X^B_{\bar{\mathcal{K}}}\cdot \bar{\varepsilon}, $
and $\tau_{\bar{\mathcal{K}}}\cdot T\pi_{/G}\cdot T\lambda=T\bar{\lambda}, $ since
$\textmd{Im}(T\bar{\gamma})\subset \bar{\mathcal{K}}. $
From the nonholonomic reduced distributional magnetic Hamiltonian equation (5.3),
$\mathbf{i}_{X^B_{\bar{\mathcal{K}}}}\omega^B_{\bar{\mathcal{K}}} =
\mathbf{d}h_{\bar{\mathcal{K}}},$ we have that $X^B_{\bar{\mathcal{K}}}
=\tau_{\bar{\mathcal{K}}}\cdot X^B_{h_{\bar{\mathcal{K}}}},$
where $ X^B_{h_{\bar{\mathcal{K}}}}$ is the magnetic Hamiltonian vector field of
the function $h_{\bar{\mathcal{K}}}: \bar{M}(\subset T^* Q/G)\rightarrow \mathbb{R}.$
Note that $\varepsilon: T^* Q \rightarrow T^* Q $ is
symplectic with respect to $\omega^B$, and
$\bar{\varepsilon}=\pi_{/G}(\varepsilon): T^* Q \rightarrow T^* Q/G$
is also symplectic along $\bar{\varepsilon}$, and
hence $X^B_{h_{\bar{\mathcal{K}}}}\cdot \bar{\varepsilon}
= T\bar{\varepsilon} \cdot X^B_{h_{\bar{\mathcal{K}}} \cdot
\bar{\varepsilon}}, $ along $\bar{\varepsilon}$, and hence
$X^B_{\bar{\mathcal{K}}}\cdot \bar{\varepsilon}
=\tau_{\bar{\mathcal{K}}}\cdot X^B_{h_{\bar{\mathcal{K}}}} \cdot\bar{\varepsilon}
= \tau_{\bar{\mathcal{K}}}\cdot T\bar{\varepsilon} \cdot X^B_{h_{\bar{\mathcal{K}}}
\cdot \bar{\varepsilon}}, $ along $\bar{\varepsilon}$.
Then we have that
\begin{align*}
& \omega^B_{\bar{\mathcal{K}}}(T\bar{\gamma} \cdot X^\varepsilon, \;
\tau_{\bar{\mathcal{K}}}\cdot T\pi_{/G}\cdot w)-
\omega^B_{\bar{\mathcal{K}}}( X^B_{\bar{\mathcal{K}}}\cdot \bar{\varepsilon},
\; \tau_{\bar{\mathcal{K}}}\cdot T\pi_{/G} \cdot w) \nonumber \\
& = -\omega^B_{\bar{\mathcal{K}}}( X^B_{\bar{\mathcal{K}}} \cdot
\bar{\varepsilon}, \; T\bar{\lambda} \cdot w)
+ \omega^B_{\bar{\mathcal{K}}}(T\bar{\lambda}\cdot X^B_H \cdot \varepsilon,
\; T\bar{\lambda} \cdot w)\\
& = \omega^B_{\bar{\mathcal{K}}}(T\bar{\lambda}\cdot X^B_H \cdot \varepsilon
- \tau_{\bar{\mathcal{K}}}\cdot T\bar{\varepsilon} \cdot
X^B_{h_{\bar{\mathcal{K}}}\cdot \bar{\varepsilon}},
\; T\bar{\lambda} \cdot w).
\end{align*}
Because the nonholonomic reduced distributional two-form
$\omega^B_{\bar{\mathcal{K}}}$ is non-degenerate, it follows that the equation
$T\bar{\gamma}\cdot X^\varepsilon = X^B_{\bar{\mathcal{K}}}\cdot
\bar{\varepsilon},$ is equivalent to the equation $T\bar{\lambda}\cdot X^B_H \cdot \varepsilon
= \tau_{\bar{\mathcal{K}}}\cdot T\bar{\varepsilon} \cdot X^B_{h_{\bar{\mathcal{K}}}
\cdot \bar{\varepsilon}}. $
Thus, $\varepsilon$ and $\bar{\varepsilon}$ satisfy the equation
$T\bar{\lambda}\cdot X^B_H \cdot \varepsilon
= \tau_{\bar{\mathcal{K}}}\cdot T\bar{\varepsilon} \cdot X^B_{h_{\bar{\mathcal{K}}}
\cdot \bar{\varepsilon}}, $ if and only if they satisfy
the Type II of Hamilton-Jacobi equation
$T\bar{\gamma}\cdot X^\varepsilon = X^B_{\bar{\mathcal{K}}}\cdot
\bar{\varepsilon} .$
\hskip 0.3cm $\blacksquare$\\

For a given nonholonomic reducible magnetic Hamiltonian system
$(T^*Q,G,\omega^B,\mathcal{D},H)$ with the associated
distributional magnetic Hamiltonian system
$(\mathcal{K},\omega^B_{\mathcal {K}},H_{\mathcal{K}})$
and the nonholonomic reduced distributional magnetic Hamiltonian system
$(\bar{\mathcal{K}},\omega^B_{\bar{\mathcal{K}}},h_{\bar{\mathcal{K}}})$,
we know that the nonholonomic dynamical vector field
$X^B_{\mathcal{K}}$ and the nonholonomic reduced dynamical
vector field $X^B_{\bar{\mathcal{K}}}$ are $\pi_{/G}$-related,
that is, $X^B_{\bar{\mathcal{K}}}\cdot \pi_{/G}=T\pi_{/G}\cdot
X^B_{\mathcal{K}}. $ Then we can prove the following Theorem 5.4,
which states the relationship between the solutions of Type II of
Hamilton-Jacobi equations and nonholonomic reduction.

\begin{theo}
For a given nonholonomic reducible magnetic Hamiltonian system
$(T^*Q,G,\omega^B,\mathcal{D},H)$ with the associated
distributional magnetic Hamiltonian system
$(\mathcal{K},\omega^B_{\mathcal {K}},H_{\mathcal{K}})$
and the nonholonomic reduced distributional magnetic Hamiltonian system
$(\bar{\mathcal{K}},\omega^B_{\bar{\mathcal{K}}},h_{\bar{\mathcal{K}}})$, assume that
$\gamma: Q \rightarrow T^*Q$ is an one-form on $Q$, and
$\lambda=\gamma \cdot \pi_{Q}: T^* Q \rightarrow T^* Q, $
and $\varepsilon: T^* Q \rightarrow T^* Q $ is a $G$-invariant symplectic map
with respect to $\omega^B$.
Moreover, assume that $\textmd{Im}(\gamma)\subset
\mathcal{M}, $ and it is $G$-invariant, $\varepsilon(\mathcal{M})\subset \mathcal{M}$,
$\textmd{Im}(T\gamma)\subset \mathcal{K}, $ and
$\bar{\gamma}=\pi_{/G}(\gamma): Q \rightarrow T^* Q/G $, and
$\bar{\lambda}=\pi_{/G}(\lambda): T^* Q \rightarrow T^* Q/G, $ and
$\bar{\varepsilon}=\pi_{/G}(\varepsilon): T^* Q \rightarrow T^* Q/G. $ Then $\varepsilon$
is a solution of the Type II of Hamilton-Jacobi equation, $T\gamma\cdot
X^\varepsilon= X^B_{\mathcal{K}}\cdot \varepsilon, $ for the distributional
magnetic Hamiltonian system $(\mathcal{K},\omega^B_{\mathcal{K}},H_{\mathcal{K}})$, if and
only if $\varepsilon$ and $\bar{\varepsilon}$ satisfy the Type II of
Hamilton-Jacobi equation $T\bar{\gamma}\cdot X^\varepsilon =
X^B_{\bar{\mathcal{K}}}\cdot \bar{\varepsilon}, $ for the nonholonomic reduced
distributional magnetic Hamiltonian system $ (\bar{\mathcal{K}},
\omega^B_{\bar{\mathcal{K}}}, h_{\bar{\mathcal{K}}} ). $
\end{theo}

\noindent{\bf Proof: } Note that
$\textmd{Im}(\gamma)\subset \mathcal{M},$ and
it is $G$-invariant, $\textmd{Im}(T\gamma)\subset \mathcal{K}, $
and hence $\textmd{Im}(T\bar{\gamma})\subset \bar{\mathcal{K}}, $ in
this case, $\pi^*_{/G}\cdot\omega^B_{\bar{\mathcal{K}}}\cdot \tau_{\bar{\mathcal{K}}}= \tau_{\mathcal{U}}\cdot
\omega^B_{\mathcal{M}}= \tau_{\mathcal{U}}\cdot i_{\mathcal{M}}^*\cdot
\omega^B, $ along $\textmd{Im}(T\bar{\gamma})$, and
$\tau_{\bar{\mathcal{K}}}\cdot T\bar{\gamma}= T\bar{\gamma},
\; \tau_{\bar{\mathcal{K}}}\cdot X^B_{\bar{\mathcal{K}}}= X^B_{\bar{\mathcal{K}}}. $
Since nonholonomic vector field $X^B_{\mathcal{K}}$ and the
vector field $X^B_{\bar{\mathcal{K}}}$ are $\pi_{/G}$-related,
that is, $X^B_{\bar{\mathcal{K}}}\cdot \pi_{/G}=T\pi_{/G}\cdot
X^B_{\mathcal{K}}, $ using the non-degenerate and nonholonomic reduced
distributional two-form $\omega^B_{\bar{\mathcal{K}}}$,
we have that
\begin{align*}
& \omega^B_{\bar{\mathcal{K}}}(T\bar{\gamma} \cdot X^\varepsilon
- X^B_{\bar{\mathcal{K}}}\cdot \bar{\varepsilon}, \; \tau_{\bar{\mathcal{K}}}\cdot T\pi_{/G}\cdot w)\\
& = \omega^B_{\bar{\mathcal{K}}}(T\bar{\gamma} \cdot X^\varepsilon, \;
\tau_{\bar{\mathcal{K}}}\cdot T\pi_{/G}\cdot w)-
\omega^B_{\bar{\mathcal{K}}}(X^B_{\bar{\mathcal{K}}}\cdot \bar{\varepsilon},
\; \tau_{\bar{\mathcal{K}}}\cdot T\pi_{/G} \cdot w) \\
& = \omega^B_{\bar{\mathcal{K}}}(\tau_{\bar{\mathcal{K}}}\cdot T\bar{\gamma}\cdot X^
\varepsilon, \; \tau_{\bar{\mathcal{K}}}\cdot T\pi_{/G}\cdot w)
-\omega^B_{\bar{\mathcal{K}}}(\tau_{\bar{\mathcal{K}}}\cdot X^B_{\bar{\mathcal{K}}}
\cdot \pi_{/G}\cdot \varepsilon, \; \tau_{\bar{\mathcal{K}}}\cdot T\pi_{/G}\cdot w)\\
& = \omega^B_{\bar{\mathcal{K}}}\cdot \tau_{\bar{\mathcal{K}}}(T\pi_{/G}\cdot T\gamma \cdot X^
\varepsilon, \; T\pi_{/G} \cdot w)
- \omega^B_{\bar{\mathcal{K}}}\cdot \tau_{\bar{\mathcal{K}}}(T\pi_{/G}\cdot
X^B_{\mathcal{K}}\cdot \varepsilon, \; T\pi_{/G}\cdot w)\\
& = \pi^*_{/G}\cdot\omega^B_{\bar{\mathcal{K}}}\cdot \tau_{\bar{\mathcal{K}}}(T\gamma \cdot X^
\varepsilon, \; w)
- \pi^*_{/G}\cdot\omega^B_{\bar{\mathcal{K}}}\cdot \tau_{\bar{\mathcal{K}}}(X^B_{\mathcal{K}} \cdot
\varepsilon, \; w)\\
& = \tau_{\mathcal{U}}\cdot i_{\mathcal{M}}^* \cdot
\omega^B (T\gamma \cdot X^
\varepsilon, \; w)- \tau_{\mathcal{U}}\cdot i_{\mathcal{M}}^* \cdot
\omega^B (X^B_{\mathcal{K}} \cdot \varepsilon, \; w).
\end{align*}
In the case we considered that $\tau_{\mathcal{U}}\cdot i_{\mathcal{M}}^* \cdot
\omega^B=\tau_{\mathcal{K}}\cdot i_{\mathcal{M}}^* \cdot
\omega^B= \omega^B_{\mathcal{K}}\cdot \tau_{\mathcal{K}}, $
and
$\tau_{\mathcal{K}}\cdot T\gamma =T\gamma, \; \tau_{\mathcal{K}} \cdot X^B_{\mathcal{K}}
= X^B_{\mathcal{K}}$,
since $\textmd{Im}(\gamma)\subset
\mathcal{M}, $ and $\textmd{Im}(T\gamma)\subset \mathcal{K}. $
Thus, we have that
\begin{align*}
& \omega^B_{\bar{\mathcal{K}}}(T\bar{\gamma} \cdot X^\varepsilon
- X^B_{\bar{\mathcal{K}}}\cdot \bar{\varepsilon}, \; \tau_{\bar{\mathcal{K}}}\cdot T\pi_{/G}\cdot w)\\
& = \omega^B_{\mathcal{K}}\cdot \tau_{\mathcal{K}}(T\gamma \cdot X^
\varepsilon, \; w)- \omega^B_{\mathcal{K}}\cdot \tau_{\mathcal{K}}(X^B_{\mathcal{K}} \cdot \varepsilon, \; w)\\
& = \omega^B_{\mathcal{K}}(\tau_{\mathcal{K}} \cdot T\gamma \cdot X^
\varepsilon, \; \tau_{\mathcal{K}} \cdot w)
- \omega^B_{\mathcal{K}}(\tau_{\mathcal{K}} \cdot X^B_{\mathcal{K}}\cdot \varepsilon,
 \; \tau_{\mathcal{K}} \cdot w)\\
& = \omega^B_{\mathcal{K}}(T\gamma \cdot X^
\varepsilon- X^B_{\mathcal{K}}\cdot \varepsilon, \; \tau_{\mathcal{K}} \cdot w).
\end{align*}
Because the distributional two-form $\omega^B_{\mathcal{K}}$
and the nonholonomic reduced distributional
two-form $\omega^B_{\bar{\mathcal{K}}}$ are both non-degenerate,
it follows that the equation
$T\bar{\gamma}\cdot X^\varepsilon=
X^B_{\bar{\mathcal{K}}}\cdot \bar{\varepsilon}, $
is equivalent to the equation $T\gamma\cdot X^\varepsilon= X^B_{\mathcal{K}}\cdot \varepsilon. $
Thus, $\varepsilon$ is a solution of the Type II of Hamilton-Jacobi equation
$T\gamma\cdot X^\varepsilon= X^B_{\mathcal{K}}\cdot \varepsilon, $ for the distributional
magnetic Hamiltonian system $(\mathcal{K},\omega^B_{\mathcal {K}},H_{\mathcal{K}})$, if and only if
$\varepsilon$ and $\bar{\varepsilon} $ satisfy the Type II of Hamilton-Jacobi
equation $T\bar{\gamma}\cdot X^\varepsilon=
X^B_{\bar{\mathcal{K}}}\cdot \bar{\varepsilon}, $ for the
nonholonomic reduced distributional magnetic Hamiltonian system
$(\bar{\mathcal{K}},\omega^B_{\bar{\mathcal{K}}},h_{\bar{\mathcal{K}}})$.
\hskip 0.3cm
$\blacksquare$

\begin{rema}
It is worthy of noting that,
the Type I of Hamilton-Jacobi equation
$T\bar{\gamma}\cdot X^ \gamma = X^B_{\bar{\mathcal{K}}}\cdot
\bar{\gamma}, $ is the equation of
the nonholonomic reduced differential one-form $\bar{\gamma}$; and
the Type II of Hamilton-Jacobi equation $T\bar{\gamma}\cdot X^\varepsilon
= X^B_{\bar{\mathcal{K}}}\cdot \bar{\varepsilon},$ is the equation of the symplectic
diffeomorphism map $\varepsilon$ and the nonholonomic reduced symplectic
diffeomorphism map $\bar{\varepsilon}$.
When $B=0$,
in this case the magnetic symplectic form $\omega^B$
is just the canonical symplectic form $\omega$ on $T^*Q$, and
the nonholonomic reducible magnetic Hamiltonian system
is just the nonholonomic reducible Hamiltonian system itself, and
the nonholonomic reduced distributional magnetic Hamiltonian system
is just the nonholonomic reduced distributional Hamiltonian system.
From the above Type I and Type II of Hamilton-Jacobi theorems, that is,
Theorem 5.2 and Theorem 5.3, we can get the Theorem 4.2 and Theorem 4.3
given in Le\'{o}n and Wang \cite{lewa15}.
It shows that Theorem 5.2 and Theorem 5.3 can be regarded as an extension of two types of
Hamilton-Jacobi theorem for the nonholonomic reduced distributional Hamiltonian system to that for the
nonholonomic reduced distributional magnetic Hamiltonian system.
\end{rema}

In order to describe the impact of different geometric structures
and constraints for the dynamics of a Hamiltonian system,
in this paper, we study the Hamilton-Jacobi theory for
the magnetic Hamiltonian system, the nonholonomic
magnetic Hamiltonian system and
the nonholonomic reducible magnetic Hamiltonian system
on a cotangent bundle, by using the distributional magnetic Hamiltonian system
and the nonholonomic reduced distributional magnetic Hamiltonian system,
which are the development of the Hamilton-Jacobi theory for
the nonholonomic Hamiltonian system and the nonholonomic reducible
Hamiltonian system given in Le\'{o}n and Wang \cite{lewa15}.
These research works reveal from the geometrical point of view the internal relationships of
the magnetic symplectic form, nonholonomic constraint,
non-degenerate distributional two form and dynamical vector fields
of a nonholonomic magnetic Hamiltonian system and the nonholonomic
reducible magnetic Hamiltonian system.
It is worthy of noting that, Marsden et al. in \cite{mawazh10} set up the
regular reduction theory of regular controlled Hamiltonian systems on a symplectic fiber
bundle, by using momentum map and the associated reduced symplectic
forms, and from the viewpoint of completeness of Marsden-Weinstein symplectic
reduction, and some developments around the above work are given in Wang and
Zhang \cite{wazh12}, Ratiu and Wang \cite{rawa12}, and Wang \cite{wa18, wa13,
wa15a, wa17}.
Since the Hamilton-Jacobi theory
is developed based on the Hamiltonian picture of dynamics, it is
natural idea to extend the Hamilton-Jacobi theory to the (regular)
controlled (magnetic) Hamiltonian systems and their a variety of reduced systems,
and it is also possible to describe the relationship between the
RCH-equivalence for controlled Hamiltonian systems and the solutions
of corresponding Hamilton-Jacobi equations, see Wang \cite{wa13d,
wa20a, wa13e} for more details.
Thus, our next topic is how to set up and develop the
nonholonomic reduction and Hamilton-Jacobi theory for the nonholonomic
controlled (magnetic) Hamiltonian systems and the distributional controlled
(magnetic) Hamiltonian systems, by analyzing carefully the geometrical
and topological structures of the phase spaces of these systems.
It is the key thought of the researches of geometrical mechanics
of the professor Jerrold E. Marsden to explore and reveal the deeply internal
relationship between the geometrical structure of phase space and the dynamical
vector field of a mechanical system. It is also our goal of pursuing and inheriting.
In addition, we note also that there have been a lot of
beautiful results of reduction theory of Hamiltonian systems in
celestial mechanics, hydrodynamics and plasma physics. Thus, it is
an important topic to study the application of reduction theory
and Hamilton-Jacobi theory of the systems in celestial mechanics, hydrodynamics
and plasma physics. These are our goals in future research.\\

\noindent {\bf Acknowledgments:}
The year of 2021 is S.S. Chern's year of Nankai University in China,
for the 110th anniversary of the birth of Professor S.S. Chern.
Shiing Shen Chern (1911-10-28---2004-12-03) is a great geometer
and a model of excellent scientists. All of the differential geometry
theory and methods for us used in the research of geometrical mechanics,
are studied from his books and his research papers.
He is a good example of us learning from him forevermore.


\begin{thebibliography}{99}

\bibitem{abma78}
Abraham R., Marsden J.E., Foundations of Mechanics, second ed.,
Addison-Wesley, Reading, MA, (1978).
\bibitem{ar89}
Arnold V.I., Mathematical Methods of Classical Mechanics, second
ed., in: Graduate Texts in Mathematics, vol. 60, Springer-Verlag,
(1989).
\bibitem{basn93}
Bates L. and $\acute{S}$niatycki J., Nonholonomic reduction, Rep.
Math. Phys. 32, 99-115(1993).
\bibitem{calemama99}
Cantrijn F., de Le\'{o}n M., Marrero J.C. and Martin de
Diego D., Reduction of constrained systems with symmetries, J. Math.
Phys., 40(2), 795-820(1999).
\bibitem{cagrmamamuro06}
Cari\~{n}ena J.F., Gr\`{a}cia X., Marmo G., Mart\'{\i}nez E.,
Mu\~{n}oz-Lecanda M. and Rom\'{a}n-Roy N., Geometric Hamilton-Jacobi
theory, Int. J. Geom. Methods Mod. Phys. 3, 1417-1458(2006).
\bibitem{cagrmamamuro10}
Cari\~{n}ena J.F., Gr\`{a}cia X., Marmo G., Mart\'{\i}nez E.,
Mu\~{n}oz-Lecanda M. and Rom\'{a}n-Roy N., Geometric Hamilton-Jacobi
theory for nonholonomic dynamical systems, Int. J. Geom. Methods
Mod. Phys. 7, 431-454(2010).
\bibitem{cemara01}
Cendra H., Marsden J.E. and Ratiu T.S., Geometric mechanics,
Lagrangian reduction and nonholonomic systems, In "Mathematics
Unlimited 2001 and Beyond" (eds. B. Engquist and W. Schmid),
Springer-Verlag, New York, 221-273(2001).
\bibitem{cudusn10}
Cushman R., Duistermaat H. and $\acute{S}$niatycki J., Geometry of
Nonholonomic Constrained Systems, Advanced series in nonlinear
dynamics, 26, (2010).
\bibitem{cukesnba95}
Cushman R., Kemppainen D., $\acute{S}$niatycki J. and Bates L.,
Geometry of nonholonomic constraints, Rep. Math. Phys.
36(2/3), 275-286(1995).
\bibitem{gema88}
Ge Z. and Marsden J.E., Lie-Poisson integrators and Lie-Poisson
Hamilton-Jacobi theory, Phys. Lett. A, 133, 134-139(1988).
\bibitem{ko92}
Koiller J., Reduction of some classical non-holonomic systems with
symmetry, Arch. Rational Mech. Anal. 118, 113-148(1992).
\bibitem{laor09}
L\'{a}zaro-Cam\'{i} J-A and Ortega J-P, The stochastic
Hamilton-Jacobi equation, J. Geom. Mech. 1, 295-315(2009).
\bibitem{lero89}
Le\'{o}n M. and Rodrigues P.R., Methods of Differential
Geometry in Analytical Mechanics, North-Holland, Amsterdam, (1989).
\bibitem{lewa15}
Le\'{o}n M. and Wang H., Hamilton-Jacobi equations for
nonholonomic reducible Hamiltonian systems on a cotangent bundle,
(arXiv: 1508.07548, a revised version).
\bibitem{lima87}
Libermann P. and Marle C.M., Symplectic Geometry and Analytical
Mechanics, Kluwer Academic Publishers, (1987).
\bibitem{ma92}
Marsden J.E., Lectures on Mechanics, in: London Mathematical Society
Lecture Notes Series, vol. 174, Cambridge University Press, (1992).
\bibitem{mamiorpera07}
Marsden J.E., Misiolek G., Ortega J.P., Perlmutter M. and Ratiu T.S.,
Hamiltonian Reduction by Stages, in: Lecture Notes in Mathematics,
vol. 1913, Springer, (2007).
\bibitem{mamora90}
Marsden J.E., Montgomery R. and Ratiu T.S., Reduction, Symmetry and
Phases in Mechanics, in: Memoirs of the American Mathematical
Society, vol. 88, American Mathematical Society, Providence, Rhode
Island, (1990).
\bibitem{mara99}
Marsden J.E. and Ratiu T.S., Introduction to Mechanics and Symmetry,
second ed., in: Texts in Applied Mathematics, vol. 17,
Springer-Verlag, New York, (1999).
\bibitem{mawazh10}
Marsden J.E., Wang H. and Zhang Z.X., Regular reduction of controlled
Hamiltonian system with symplectic structure and symmetry, Diff.
Geom. Appl., 33(3), 13-45(2014), (arXiv: 1202.3564).
\bibitem{mawe74}
Marsden J.E. and Weinstein A., Reduction of symplectic manifolds with
symmetry, Rep. Math. Phys. 5, 121-130(1974).
\bibitem{orra04}
Ortega J-P and Ratiu T.S., Momentum Maps and Hamiltonian Reduction,
in: Progress in Mathematics, vol. 222, Birkh\"{a}user, (2004).
\bibitem{pa07}
Patrick G.W., Variational development of the semi-symplectic
geometry of nonholonomic mechanics, Rep. Math. Phys. 59,
145-184(2007).
\bibitem{rawa12}
Ratiu T.S. and Wang H., Poisson reduction by controllability distribution
for a controlled Hamiltonian system, (arXiv: 1312.7047).
\bibitem{wa18}
Wang H. Some developments of reduction theory for
controlled Hamiltonian system with symmetry ( in Chinese),
Sci Sin Math, 2018, \text{48}, 1-12.
\bibitem{wa13}
Wang H. Reductions of controlled Hamiltonian system with symmetry,
In: Symmetries and Groups in Contemporary Physics,
( Bai C M. Gazeau J P. and Ge M L. eds),
World Scientific, 2013, 639-642.
\bibitem{wa15a}
Wang H., Regular reduction of a controlled magnetic
Hamiltonian system with symmetry of the Heisenberg group,
(arXiv: 1506.03640, a revised version).
\bibitem{wa17}
Wang H., Hamilton-Jacobi theorems for regular reducible Hamiltonian
systems on a cotangent bundle, Jour. Geom. Phys., \textbf{119}
82-102, (2017).
\bibitem{wa13d}
Wang H., Hamilton-Jacobi equations for a regular controlled Hamiltonian
system and its reduced systems, (arXiv: 1305.3457, a revised version),
To appear in Acta Mathematica Scientia, English Series, 2022.
\bibitem{wa20a}
Wang H., Dynamical equations of the controlled rigid spacecraft with a rotor,
(arXiv: 2005.02221).
\bibitem{wa13e}
Wang H., Symmetric reduction and Hamilton-Jacobi equations for the controlled
underwater vehicle-rotor system, ( arXiv: 1310.3014, a revised version ).
\bibitem{wazh12}
Wang H. and Zhang Z.X., Optimal reduction of controlled Hamiltonian
system with Poisson structure and symmetry, Jour. Geom. Phys., 62
(5), 953-975(2012).
\bibitem{wo92}
Woodhouse N.M.J., Geometric Quantization, second ed., Clarendon
Press, Oxford, (1992).

\end{thebibliography}
\end{document}